\documentclass[12pt]{amsart}

\usepackage[margin=1.2in]{geometry}
\usepackage[english]{babel}
\usepackage[utf8]{inputenc}
\usepackage{fancyhdr}
\usepackage{natbib}
\usepackage{amsmath}

\RequirePackage{amsthm,amsmath,amsfonts,amssymb}
\usepackage{graphicx}
\usepackage{hyperref}
\usepackage{graphicx}
\usepackage{placeins}
\usepackage{caption}
\usepackage{subcaption}
\usepackage{booktabs}
\numberwithin{equation}{section}

\usepackage{setspace}

\theoremstyle{plain}
\theoremstyle{remark} 
\numberwithin{equation}{section} 

\title{A unified study for estimation of order restricted location/scale parameters under the generalized Pitman nearness criterion}
\author{Naresh Garg  and Neeraj Misra \\ {\footnotesize Department of Mathematics and Statistics\\Indian Institute of Technology Kanpur \\Kanpur-208016, Uttar Pradesh, India}}

\makeatletter
\def\@seccntformat#1{%
  \protect\textup{\protect\@secnumfont
    \ifnum\pdfstrcmp{subsection}{#1}=0 \bfseries\fi% subsection # in \bfseries
    \csname the#1\endcsname
    \protect\@secnumpunct
  }%
}  
\makeatother

\begin{document}
\maketitle
\section*{\textbf{Abstract}}

We consider component-wise estimation of order restricted location/scale parameters of a general bivariate location/scale distribution under the generalized Pitman nearness criterion (GPN). We develop some general results that, in many situations, are useful in finding improvements over location/scale equivariant estimators. In particular, under certain conditions, these general results provide improvements over the unrestricted Pitman nearest location/scale equivariant estimators and restricted maximum likelihood estimators. The usefulness of the obtained results is illustrated through their applications to specific probability models. A simulation study has been considered to compare how well different estimators perform under the GPN criterion with a specific loss function. 
\\~\\ \textbf{Keywords:} Generalised Pitman Nearness (GPN) Criterion; Location equivariant estimator; Pitman Nearness (PN) Criterion; Restricted Parameter Space; Scale equivariant estimator; Unrestricted Parameter Space.

%%\pacs[JEL Classification]{D8, H51}

%%\pacs[MSC Classification]{35A01, 65L10, 65L12, 65L20, 65L70}

	\section{\textbf{Introduction}}  \label{introduction}	
The problem of estimating real-valued parameters of a set of distributions, when it is known apriori that these parameters follow certain order restrictions, is of great relevance. 
For example, in a clinical trial, where estimating average blood pressures of two groups of Hypertension patients, one treated with a standard drug and the other with a placebo, is of interest, it can be assumed that the average blood pressure of Hypertension patients treated with the standard drug is lower than the average blood pressure of Hypertension patients treated with the placebo. 
%Taketomi (\citeyear{taketomi2021meta}) provides a compelling example of how estimation problems for order-restricted parameters are applied in meta-analysis.
Such estimation problems have been extensively studied in the literature. For a detailed account of work carried out in this area, one may refer to Barlow et al. (\citeyear{MR0326887}), Robertson et al. (\citeyear{MR961262}) and van Eeden (\citeyear{MR2265239}). \vspace*{2mm}

Early work on estimation of order restricted parameters deals with obtaining isotonic regression estimators or restricted maximum likelihood estimators (MLE) (see Brunk (\citeyear{MR73894}), van Eeden (\citeyear{MR0083859}, \citeyear{MR0090197}, \citeyear{MR0090196}, \citeyear{MR0102874}) and Robertson et al. (\citeyear{MR961262})).
%After that, researchers found that the invariance approach to statistical inference is also efficient in deriving better estimators. However, in most of the studies, specific distributions with independent marginals and specific loss functions have been considered. 
%Also, most studies in this area deal with shrinking the risk of location/scale equivariant estimators to a minimum risk and yields estimators that improve over the unrestricted BLEE/BSEE and/or unrestricted MLE. 
Subsequently, a lot of work was carried out using the decision theoretic approach under different loss functions. Some of the key contributions in this direction are due to  Katz (\citeyear{MR150874}), Cohen and Sackrowitz (\citeyear{MR270483}), Brewster-Zidek (\citeyear{MR381098}), Lee (\citeyear{MR615447}), Kumar and Sharma (\citeyear{MR981031}, \citeyear{MR1058934},\citeyear{MR1165706}), Kelly (\citeyear{MR994278}), Kushary and Cohen (\citeyear{MR1029476}), Kaur and Singh (\citeyear{MR1128873}), Gupta and Singh (\citeyear{gupta1992}), Vijayasree and Singh (\citeyear{MR1220404}), Hwang and Peddada (\citeyear{MR1272076}), Kubokawa and Saleh (\citeyear{MR1370413}), Vijayasree et al. (\citeyear{MR1345425}), Misra and Dhariyal (\citeyear{MR1326266}), Garren (\citeyear{MR1802627}), Misra et al. (\citeyear{MR1904424}, \citeyear{MR2205815}), Peddada et al. (\citeyear{MR2202656}), Chang and Shinozaki (\citeyear{chang2015}) and Patra (\citeyear{Patra}).

% For a detailed account of other developments in this direction, one may refer to Van Eeden (\citeyear{MR2265239}).

% Also, in most of the studies, a minimizing risk function criterion is used to yield an improved estimator that improve over the unrestricted BLEE/BSEE and/or unrestricted MLE. In this paper, we use the generalised Pitman nearness (GPN) criterion to derive improved estimators under the restricted parameter space and further compare these estimators with the usual estimators (BLEE or MLE) in terms of the GPN criterion. Best of our knowledge, the GPN criterion is not yet used in the literature to compare and obtain improved estimators for order-restricted parameters. \vspace*{2mm}	

A popular alternative criterion to compare different estimators is the Pitman nearness criterion, due to Pitman (\citeyear{pitman1937}). This criterion of comparing two estimators is based on the probability that one estimator is closer to the estimand than the other estimator under the $L_1$ distance (i.e., absolute error loss function). Rao (\citeyear{rao1981}) pointed out advantages of the Pitman nearness (PN) criterion over the mean squared error criterion. Keating (\citeyear{keating1985}) further advocated Rao's findings through certain examples. Keating and Mason (\citeyear{keating1985m}) provided some practical examples where the PN criterion is more relevant than minimizing the risk function. Also, Peddada (\citeyear{peddada1985}) and Rao et al. (\citeyear{MR860477}) extended the notion of PN criterion by defining the generalised Pitman Criterion (GPN) based on a general loss function (in place of the absolute error loss function). For a detailed study of the PN criterion, one may check out the monograph by Keating et al. (\citeyear{keating1993}). \vspace*{2mm}

The PN criterion has been extensively used in the literature for different estimation problems (see Nayak (\citeyear{nayak1990}) and Keating (1993)). However, there are only a limited number of studies on use of the PN criterion in estimating order restricted parameters. For component-wise estimation of order restricted means of two independent normal distributions, having a common unknown variance, Gupta and Singh (\citeyear{gupta1992}) showed that restricted MLEs are nearer to the respective population means than unrestricted MLEs under the PN criterion. Analogous result was also proved for the estimation of the common variance, which is also considered in Misra et al. (\citeyear{MR2205815}). Chang et al. (\citeyear{chang2017}, \citeyear{chang2020}) considered estimation of order restricted means of a bivariate normal distribution, having a known covariance matrix, and established that, under a modified PN criterion, the restricted MLEs are nearer to respective population means than some of the estimators proposed by Hwang and Peddada (\citeyear{MR1272076}) and Tan and Peddada (\citeyear{tan2000}). Ma and Liu (\citeyear{ma2014}) considered the problem of estimating order restricted scale parameters of two independent gamma distributions. Under the PN criterion, they have compared restricted MLEs of scale parameters with the best unbiased estimators. Some other studies in this direction are due to Misra and van der Meulen (\citeyear{misra1997}), Misra et al. (\citeyear{MR2205815}), and Chang and Shinozaki (\citeyear{chang2015}). Most of these studies are centered around specific probability distributions (mostly, normal and gamma) and the absolute error loss in the PN criterion. In this paper, we aim to unify these studies by considering the problem of estimating order restricted location/scale parameters of a general bivariate location/scale model under the GPN criterion with a general loss function. We will develop some general results that, in certain situations, are useful in finding improvements over location/scale equivariant estimators. As a consequence of these general results, we will obtain estimators improving upon the unrestricted Pitman nearest location/scale equivariant estimators (PNLEE/PNSEE). We also consider some applications of these general results to specific probability models and obtain improvements over the PNLEE/PNSEE and the restricted MLEs. 

% Gupta and Singh (\citeyear{gupta1992}) showed that the restricted MLE is nearer to the parameter than the BLEE under the PN criterion for estimating the ordered means of two independent normal distributions. Further, Chang et al. (\citeyear{chang2017},\citeyear{chang2020}) showed that, under the modified PN criterion, the restricted MLE is nearer than the Hwang and Peddada (\citeyear{MR1272076}) estimator and the Peddada et al. (\citeyear{MR2202656}) estimator for estimating ordered means of a bivariate normal distribution with a known covariance matrix. Ma and Liu (\citeyear{ma2014}) consider the problem of estimating the order restricted scale parameters of two independent gamma distributions. Under the PN criterion, they have compared the restricted MLEs of scale parameters with unbiased estimators (UEs). Some other studies in this area are due to Misra and Meulen (\citeyear{misra1997}), Misra and Singh (\citeyear{misra2004}), and Chang and Shinozaki (\citeyear{chang2015}). All the above studies used the PN criterion only for comparison of estimators of order restricted parameters. In the literature, it seems that no work has been done regarding using the PN criterion to obtain improved estimators when parameters are order restricted. At the same time, many authors used this criterion with an invariance approach to derive improved estimators when parameter space is not restricted.
\vspace*{2mm}

%This paper considers estimating order restricted location/scale parameters of the general bivariate location/scale model. We use the GPN criterion to derive better estimators over any arbitrary equivariant estimator. When parameter space is unrestricted, under the GPN criterion, we obtain the Pitman nearest equivariant estimator (PNEE), and under minimizing the risk function criterion, we obtain the BLEE/BSEE. Further,  we use the GPN criterion related to the PNEE and derive an estimator which dominates the PNEE in the sense of the GN criterion under the restricted parameter space. Similarly,  we use the GPN criterion related to the BLEE/BSEE and derive an estimator which dominates BLEE/BSEE in the sense of the GN criterion under the restricted parameter space. Through the applications of our findings, we will also generalize the results of Chang et al. (\citeyear{chang2017},\citeyear{chang2020}) and Ma and Liu (\citeyear{ma2014}) under the PN criterion to the GPN criterion with a general loss function.\vspace*{2mm}

The rest of the paper is organized as follows. In Section \ref{1}, we introduce some useful notations, definitions, and results that are used later in the paper. Sections \ref{3.1} and \ref{3.2} (\ref{4.1} and \ref{4.2}), respectively, deal with estimating the smaller and larger location (scale) parameters under the criterion of GPN. In Section \ref{5}, we present a simulation study to compare performances of various competing estimators.

\section{\textbf{Some Useful Notations, Definitions and Results}} \label{1}
\noindent
The following notion of the Pitman nearness criterion was first introduced by Pitman (\citeyear{pitman1937}). 
\\~\\\textbf{Definition 2.1} Let $\bold{X}$ be a random vector having a distribution involving an unknown parameter $\boldsymbol{\theta}\in \Theta$ ($\boldsymbol{\theta}$ may be vector valued).  Let $\delta_1$ and $\delta_2$ be two estimators of a real-valued estimand $\tau(\boldsymbol{\theta})$. Then, the Pitman nearness (PN) of $\delta_1$ relative to $\delta_2$ is defined by
$$PN(\delta_1,\delta_2;\boldsymbol{\theta})=P_{\boldsymbol{\theta}}[\vert \delta_1-\tau(\boldsymbol{\theta})\vert <\vert \delta_2-\tau(\boldsymbol{\theta})\vert], \; \;\boldsymbol{\theta}\in\Theta,$$
and the estimator $\delta_1$ is said to be nearer to $\tau(\boldsymbol{\theta})$ than $\delta_2$ if $PN(\delta_1,\delta_2;\boldsymbol{\theta})\geq \frac{1}{2},\;\forall\; \boldsymbol{\theta}\in\Theta$, with strict inequality for some $\boldsymbol{\theta}\in\Theta$.
\vspace*{2mm}

Two drawbacks of the above criterion are that it does not take into account that estimators $\delta_1$ and $\delta_2$ may coincide over a subset of the sample space and that it is only based on the $L_1$ distance (absolute error loss). To take care of these deficiencies, Nayak (\citeyear{nayak1990}) and Kubokawa (\citeyear{kubokawa1991}) modified the Pitman (\citeyear{pitman1937}) nearness criterion and defined the generalized Pitman nearness (GPN) criterion based on general loss function $L(\boldsymbol{\theta},\delta).$
\\~\\\textbf{Definition 2.2}  Let $\bold{X}$ be a random vector having a distribution involving an unknown parameter $\boldsymbol{\theta}\in \Theta$ and let $\tau(\boldsymbol{\theta})$ be a real-valued estimand. Let $\delta_1$ and $\delta_2$ be two estimators of the estimand $\tau(\boldsymbol{\theta})$. Also, let $L(\boldsymbol{\theta},a)$ be a specified loss function for estimating $\tau(\boldsymbol{\theta})$. Then, the generalized Pitman nearness (GPN) of $\delta_1$ relative to $\delta_2$ is defined by
$$GPN(\delta_1,\delta_2;\boldsymbol{\theta})=P_{\boldsymbol{\theta}}[L(\boldsymbol{\theta},\delta_1) <L(\boldsymbol{\theta},\delta_2)]+\frac{1}{2} P_{\boldsymbol{\theta}}[L(\boldsymbol{\theta},\delta_1) =L(\boldsymbol{\theta},\delta_2)], \; \; \boldsymbol{\theta}\in\Theta.$$
The estimator $\delta_1$ is said to be nearer to $\tau(\boldsymbol{\theta})$ than $\delta_2$, under the GPN criterion, if $GPN(\delta_1,\delta_2;\boldsymbol{\theta})\geq \frac{1}{2},\;\forall\; \boldsymbol{\theta}\in\Theta$, with strict inequality for some $\boldsymbol{\theta}\in\Theta$.
\\~\\\textbf{Definition 2.3} Let $\mathcal{D}$ be a class of estimators of a real-valued estimand $\tau(\boldsymbol{\theta})$. Let $L(\boldsymbol{\theta},a)$ be a given loss function. Then, an estimator $\delta^{*}$ is said to be the Pitman nearest within the class $\mathcal{D}$, if 
$$GPN(\delta^{*},\delta;\boldsymbol{\theta}) \geq \frac{1}{2}, \; \forall \; \delta\in \mathcal{D},\; \boldsymbol{\theta}\in\Theta,$$
with strict inequality for some $\boldsymbol{\theta}\in\Theta$.
%	The basic definition of median of a continuous random variable (r.v.) is given below.
% \\~\\\textbf{Definition 2.4} Let $X$ be a continuous r.v. with df $F(x)$. If $F(x)=P[X\leq x]=\frac{1}{2}$ have unique solution $m$, then $m$ is the median of $X$ and, if $F(x)=\frac{1}{2}$ does not have unique solution, then $m=\inf\{x\in \Re: P[X\leq x]= \frac{1}{2}\}$ is the unique median of $X$.
\vspace*{2mm}

The following result, famously known as Chebyshev's inequality, will be used in our study (see Marshall and Olkin (\citeyear{MR2363282})).
\\~\\ \textbf{Proposition 2.1} Let $S$ be random variable (r.v.) and let $k_1(\cdot)$ and $k_2(\cdot)$ be real-valued monotonic functions defined on the distributional support of the r.v. $S$. If $k_1(\cdot)$ and $k_2(\cdot)$ are monotonic functions of the same (opposite) type, then
$$E[k_1(S)k_2(S)]\geq (\leq ) E[k_1(S)] E[k_2(S)],$$
provided the above expectations exist.

\section{\textbf{Improved Estimators for Restricted Location Parameters}} \label{3}

Let $\bold{X}=(X_1,X_2)$ be a random vector with a joint probability density function (p.d.f.)
\begin{equation}\label{eq:2.1}
	f_{\boldsymbol{\theta}}(x_1,x_2)= 	f(x_1-\theta_1,x_2-\theta_2),\; \; \;(x_1,x_2)\in \Re^2, 
\end{equation} 
where $f(\cdot,\cdot) $ is a specified Lebesgue p.d.f. on $\Re^2$ and $\boldsymbol{\theta}=(\theta_1,\theta_2)\in \Theta_0=\{(t_1,t_2)\in\Re^2:t_1 \leq t_2\}$ is the vector of unknown restricted location parameters; here $\Re$ denotes the real line and $\Re^2=\Re\times \Re$. Generally, $\bold{X}=(X_1,X_2)$ would be a minimal-sufficient statistic based on a bivariate random sample or two independent random samples, as the case may be.  \vspace{2mm}

Consider estimation of the location parameter $\theta_i$ under the GPN criterion with a general loss function $L_i(\boldsymbol{\theta},a)=W(a-\theta_i),\;\boldsymbol{\theta}\in\Theta_0,\; a\in\mathcal{A},\;i=1,2,$     
where $\mathcal{A}=\Re$ and $W:\Re\rightarrow [0,\infty)$ is a specified non-negative function such that $W(0)=0$, $W(t)$ is strictly decreasing on $(-\infty,0)$ and strictly increasing on $(0,\infty)$.	Throughout this section, whenever term "general loss function" is used, it refers to a loss function having the above properties. Also, in this section, the GPN criterion is considered with a general loss function described above.\vspace{2mm}

The problem of estimating restricted location parameter $\theta_i \,(i=1,2),$ under the GPN criterion, is invariant under the group of transformations $\mathcal{G}=\{g_c:\,c\in\Re\},$ where $g_c(x_1,x_2)$ $=(x_1+c,x_2+c),\; (x_1,x_2)\in\Re^2,\;c\in\Re.$ Under the group of transformations $\mathcal{G}$, any location equivariant estimator of $\theta_i$ has the form
\begin{equation}\label{eq:2.2}
	\delta_{\psi}(\bold{X})=X_i-\psi(D),
\end{equation}
for some function $\psi:\,\Re\rightarrow \Re\,,\;i=1,2,$ where $D=X_2-X_1$. Let $f_D(t\vert \lambda)$ be the p.d.f. of r.v. $D=X_2-X_1$, where $\lambda=\theta_2-\theta_1\in[0,\infty)$. Note that the distribution of $D$ depends on $\boldsymbol{\theta}\in\Theta_0$ only through $\lambda=\theta_2-\theta_1\in[0,\infty)$. Exploiting the prior information of order restriction on parameters $\theta_1$ and $\theta_2$ ($\theta_1\leq \theta_2$), we aim to obtain estimators that are Pitman nearer to $\theta_i,\,i=1,2$. \vspace*{2mm}

The following lemma will be useful in proving the main results of this section (also see Nayak (\citeyear{nayak1990}) and Zhou and Nayak (\citeyear{zhou2012})).
\\~\\	 \textbf{Lemma 3.1} Let $Y$ be a r.v. having the Lebesgue p.d.f. and let $m_Y$ be the median of $Y$. Let $W:\Re\rightarrow [0,\infty)$ be a function such that $W(0)=0$, $W(t)$ is strictly decreasing on $(-\infty,0)$ and strictly increasing on $(0,\infty)$. Then, for $-\infty< c_1<c_2\leq m_Y$ or $-\infty<m_Y\leq c_2<c_1$,
$GPN= P[W(Y-c_2)<W(Y-c_1)]+\frac{1}{2} P[W(Y-c_2)=W(Y-c_1)]>\frac{1}{2}$.
\begin{proof}.  We have the following two cases:
	\\~\\Case 1: $-\infty<c_1<c_2\leq m_Y<\infty$
	
	In this case $Y-m_Y\leq Y-c_2<Y-c_1$ and, thus, $Y-m_Y\geq 0$ implies that $W(Y-c_2)<W(Y-c_1)$. Consequently,
	$$GPN\geq P[W(Y-c_2) <W(Y-c_1)]> P[Y\geq m_Y]= \frac{1}{2}.$$
	\\Case 2: $-\infty<m_Y\leq c_2<c_1<\infty$
	
	In this case $Y-m_Y\geq Y-c_2>Y-c_1$. Thus, $Y-m_Y\leq 0$ implies that $W(Y-c_2)<W(Y-c_1)$. Hence
	$$GPN\geq P[W(Y-c_2) <W(Y-c_1)]> P[Y\leq m_Y]= \frac{1}{2}.$$
\end{proof}

Note that, in the unrestricted case (parameter space $\Theta=\Re^2$), the problem of estimating $\theta_i,\;i=1,2,$ under the GPN criterion is invariant under the group of transformations $\mathcal{G}_0=\{g_{c_1,c_2}:\,(c_1,c_2)\in\Re^2\},$ where $g_{c_1,c_2}(x_1,x_2)$ $=(x_1+c_1,x_2+c_2),\; (x_1,x_2)\in\Re^2,\;(c_1,c_2)\in\Re^2.$ Any location equivariant estimator is of the form $\delta_{i,c}(\bold{X})=X_i-c,\; c\in \Re,\; i=1,2$.
An immediate consequence of Lemma 3.1 is that, under the unrestricted parameter space $\Theta=\Re^2$, the Pitman nearest location equivariant estimator (PNLEE) of $\theta_i$, under the GPN criterion, is $\delta_{i,PNLEE}(\bold{X})=X_i-m_{0,i}$, where $m_{0,i}$ is the median of the r.v. $Z_i=X_i-\theta_i,\;i=1,2$.\vspace*{2mm}

% Also, in the unrestricted case ($\Theta=\Re^2$), the problem of estimating the location parameter $\theta_i$, under the loss function $L_i(a,\boldsymbol{\theta})=W(a-\theta_i),\; a\in \mathcal{A}=\Re,\; \boldsymbol{\theta}\in \Re^2,\;i=1,2$, is invariant under the additive group of transformations $\mathcal{G}_0$. The unrestricted best location equivariant estimator (BLEE) of $\theta_i$ is $\delta_{i,BLEE}(\bold{X})=X_i-c_{0,i}$, where $c_{0,i}$ is the unique solution of the equation $\int_{-\infty}^{\infty} \int_{-\infty}^{\infty}\, W'(s_i-c) \,f(s_1,s_2)\,ds_1\,ds_2=0,\;i=1,2.$

\subsection{\textbf{Estimation of the Smaller Location Parameter $\theta_1$}} \label{3.1}

Let $Z_1=X_1-\theta_1$, $\lambda=\theta_2-\theta_1$ and $f_D(t\vert \lambda)$ be the p.d.f. of r.v. $D=X_2-X_1$. Let $\delta_{\xi}(\bold{X})=X_1-\xi(D)$ and $\delta_{\psi}(\bold{X})=X_1-\psi(D)$ be two location equivariant estimators of $\theta_1$, where $\xi$ and $\psi$ are real-valued functions defined on $\Re$. Then, the GPN of $\delta_{\xi}(\bold{X})$ relative to $\delta_{\psi}(\bold{X})$ is given by 
\begin{align*}
	GPN(\delta_{\xi},\delta_{\psi};\boldsymbol{\theta})&=P_{\boldsymbol{\theta}}[W(Z_1-\xi(D)) <W(Z_1-\psi(D))]
	\\&\quad +\frac{1}{2} P_{\boldsymbol{\theta}}[W(Z_1-\xi(D)) =W(Z_1-\psi(D))], \; \; \boldsymbol{\theta}\in\Theta_0,\\
	&=\int_{-\infty}^{\infty} g_{1,\lambda}(\xi(t),\psi(t),t) f_D(t\vert \lambda) dt,\;\;\lambda\geq 0,
\end{align*}
where, for $\lambda \geq 0$ and fixed $t$ in the support of the distribution of r.v. $D$, \begin{align}g_{1,\lambda}(\xi(t),\psi(t),t)&= P_{\boldsymbol{\theta}}[W(Z_1-\xi(t)) <W(Z_1-\psi(t))\vert D=t] \nonumber
	\\&\qquad +\frac{1}{2}P_{\boldsymbol{\theta}}[W(Z_1-\xi(t))\! =\!W(Z_1-\psi(t))\vert D=t].\end{align}
For any fixed $\lambda\geq 0$ and $t$, let $m_{\lambda}^{(1)}(t)$ denote the median of the conditional distribution of $Z_1$ given $D=t$. For any fixed $t$, the conditional p.d.f. of $Z_1$ given $D=t$ is $f_{\lambda}(s\vert t)=\frac{f(s,s+t-\lambda)}{f_D(t\vert \lambda)}$ and $f_D(t\vert \lambda)=\int_{-\infty}^{\infty}f(y,y+t-\lambda)dy$, $\lambda\geq 0$. Thus, $\int_{-\infty}^{m_{\lambda}^{(1)}(t)} \, f(s,s+t-\lambda)ds=\frac{1}{2} \int_{-\infty}^{\infty}f(s,s+t-\lambda)ds$. For any fixed $t$, using Lemma 3.1, we have $g_{1,\lambda}(\xi(t),\psi(t),t)> \frac{1}{2},\;\forall \;\lambda\geq 0$, if $\psi(t)<\xi(t)\leq m_{\lambda}^{(1)}(t),\;\forall \;\lambda\geq 0$, or if $ m_{\lambda}^{(1)}(t)\leq \xi(t)<\psi(t),\;\forall \;\lambda\geq 0$. Also, note that, for any fixed $t$, $g_{1,\lambda}(\psi(t),\psi(t),t)= \frac{1}{2},\;\forall \;\lambda\geq 0$. These observations, along with Lemma 3.1, yield the following result.

%Using Lemma 3.1, we get $g_{1,\lambda}(m_{\lambda}^{(1)}(t),\psi(t),t)\geq \frac{1}{2},\;\forall \;\lambda\geq 0$ and $\forall\; t$, where, for any fixed $\lambda\geq 0$ and $t$, $m_{\lambda}^{(1)}(t)$ is the median of conditional distyribution of $Z_1$ given $D=t$.

\vspace*{2mm}

%Now we develop the following theorem that, in certain situations, provides a shrinkage estimator dominating an arbitrary location equivariant estimator $\delta_{\psi}(\bold{X})=X_1-\psi(D)$.
%\vspace*{2mm}

\noindent
\textbf{Theorem 3.1.1.} Let $\delta_{\psi}(\bold{X})=X_1-\psi(D)$ be a location equivariant estimator of $\theta_1$, where $\psi:\Re\rightarrow \Re$. Let $l^{(1)}(t)$ and $u^{(1)}(t)$ be functions such that $l^{(1)}(t)\leq m_{\lambda}^{(1)}(t)\leq u^{(1)}(t), \;\forall\;\lambda\geq 0$ and any $t$. For any fixed $t$, define $\psi^{*}(t)\!=\!\max\{l^{(1)}(t),\min\{\psi(t),u^{(1)}(t)\}\}$. Then, under the GPN criterion with a general loss function,
the estimator $\delta_{\psi^{*}}(\bold{X})\!=\!X_1-\psi^{*}(D)$ is Pitman nearer to $\theta_1$ than the estimator $\delta_{\psi}(\bold{X})=X_1-\psi(D)$, for all $\boldsymbol{\theta}\in \Theta_0$, provided $P_{\boldsymbol{\theta}}[l^{(1)}(D)\leq \psi(D)\leq u^{(1)}(D)]<1,\;\forall\; \boldsymbol{\theta}\in \Theta_0$.
\begin{proof}\!\!\!.
	The GPN of the estimator $\delta_{\psi^{*}}(\bold{X})=X_1-\psi^{*}(D)$ relative to $\delta_{\psi}(\bold{X})=X_1-\psi(D)$ can be written as 
	\begin{align*}
		GPN(\delta_{\psi^{*}},\delta_{\psi};\boldsymbol{\theta})=\int_{-\infty}^{\infty} g_{1,\lambda}(\psi^{*}(t),\psi(t),t) f_D(t\vert \lambda) dt,\;\;\lambda\geq 0,
	\end{align*}
	where $g_{1,\lambda}(\cdot,\cdot,\cdot)$ is defined by (3.3).

	Let $A=\{t:\psi(t)<l^{(1)}(t)\}$, $B=\{t:l^{(1)}(t)\leq \psi(t)\leq u^{(1)}(t)\}$ and $C=\{t:\psi(t)>u^{(1)}(t)\}$. Clearly
	$$\psi^{*}(t)=\begin{cases} l^{(1)}(t), & t\in A\\  \psi(t), & t\in B\\ u^{(1)}(t), & t\in C\end{cases}.$$
	Since $l^{(1)}(t)\leq m_{\lambda}^{(1)}(t)\leq u^{(1)}(t),\; \forall\; \lambda\geq 0$ and $t$, using Lemma 3.1, we have $g_{1,\lambda}(\psi^{*}(t),\psi(t),t)> \frac{1}{2}, \; \forall \; \lambda\geq 0$, provided $t\in A\cup C$. Also, for $t\in B$, $g_{1,\lambda}(\psi^{*}(t),\psi(t),t)= \frac{1}{2}, \; \forall \; \lambda\geq 0$. Since $P_{\underline{\theta}}(A\cup C)>0,\; \forall \; \boldsymbol{\theta}\in \Theta_0$, we have
	\\~\\  $GPN(\delta_{\psi^{*}},\delta_{\psi};\boldsymbol{\theta})$
	\begin{align*}
		&\!=\!\int_{A}\! g_{1,\lambda}(\psi^{*}(t),\psi(t),t) f_D(t\vert \lambda) dt \!+ \!\int_{B}\! g_{1,\lambda}(\psi^{*}(t),\psi(t),t) f_D(t\vert \lambda) dt \!\\&\qquad\qquad\qquad\qquad\qquad\qquad\qquad\qquad\qquad\qquad\qquad+\! \int_{C}\! g_{1,\lambda}(\psi^{*}(t),\psi(t),t) f_D(t\vert \lambda) dt\\
		&> \frac{1}{2}, \;\; \boldsymbol{\theta}\in\Theta_0.
	\end{align*}
\end{proof}

Using arguments similar to those used in the proof of Theorem 3.1.1, one can, in fact, obtain a class of estimators improving over an arbitrary location equivariant estimator, under certain conditions. The proof of the following corollary is contained in the proof of Theorem 3.1.1, and hence skipped.
\\~\\	\textbf{Corollary 3.1.1.} Let $\delta_{\psi}(\bold{X})=X_1-\psi(D)$ be a location equivariant estimator of $\theta_1$ such that $P_{\boldsymbol{\theta}}[l^{(1)}(D)\leq \psi(D)\leq u^{(1)}(D)]<1,\;\forall\; \boldsymbol{\theta}\in \Theta_0$, where $l^{(1)}(\cdot)$ and $u^{(1)}(\cdot)$ are as defined in Theorem 3.1.1. Let $\psi_{1,0}:\Re\rightarrow \Re$ be such that $\psi(t)< \psi_{1,0}(t) \leq l^{(1)}(t)$, whenever $\psi(t)< l^{(1)}(t)$, and $u^{(1)}(t)\leq \psi_{1,0}(t)<\psi(t)$, whenever $u^{(1)}(t)< \psi(t)$. Also let $\psi_{1,0}(t)=\psi(t)$, whenever $l^{(1)}(t)\leq \psi(t) \leq u^{(1)}(t)$. Let $\delta_{\psi_{1,0}}(\bold{X})=X_1-\psi_{1,0}(D)$. Then, 
$GPN(\delta_{\psi_{1,0}},\delta_{\psi};\boldsymbol{\theta})> \frac{1}{2},\; \forall \; \boldsymbol{\theta}\in \Theta_0$.\vspace*{2mm}

Ironically, the result stated in Corollary 3.1.1 is general than the one stated in Theorem 3.1.1. However, among all the improved estimators provided through Corollary 3.1.1, the maximum improvement is provided by the one covered under Theorem 3.1.1.   \vspace*{2mm}

The following corollary provides improvements over the unrestricted PNLEE $\delta_{1,PNLEE}(\bold{X})=X_1-m_{0,1}$, under the restricted parameter space $\Theta_0$.
\\~\\\textbf{Corollary 3.1.2.} Let $l^{(1)}(t)$ and $u^{(1)}(t)$ be as defined in Theorem 3.1.1. Suppose that $P_{\boldsymbol{\theta}}[l^{(1)}(D)\leq m_{0,1}\leq u^{(1)}(D)]<1,\; \forall\; \boldsymbol{\theta}\in \Theta_0$. Define, for any fixed $t$, $\xi^{*}(t)\!=\!\max\{l^{(1)}(t),\min\{m_{0,1},u^{(1)}(t)\}\}$. Then, for every $\boldsymbol{\theta}\in \Theta_0$,
the estimator $\delta_{\xi^{*}}(\bold{X})\!=\!X_1-\xi^{*}(D)$ is Pitman nearer to $\theta_1$ than the PNLEE $\delta_{1,PNLEE}(\bold{X})=X_1-m_{0,1}$, under the GPN criterion.\vspace*{2mm}

In order to identify $l^{(1)}(t)$ and $u^{(1)}(t)$, satisfying $l^{(1)}(t)\leq m_{\lambda}^{(1)}(t)\leq u^{(1)}(t),\; \forall\; \lambda \geq 0 \text{ and } t$, the following lemma will be useful.           
\\~\\\textbf{Lemma 3.1.1.} If, for every fixed $\lambda\geq 0$ and $t$, $f(s,s+t-\lambda)/f(s,s+t)$ is increasing (decreasing) in $s$ (wherever the ratio is not of the form $0/0$), then, for every fixed $t$, $m_{\lambda}^{(1)}(t)$ is an increasing (decreasing) function of $\lambda\in [0,\infty)$.
\begin{proof}\!\!\!.
	Let us fix $t$, $\lambda_1$ and $\lambda_2$, such that $0\leq \lambda_1<\lambda_2<\infty$. Then, the hypothesis of the lemma implies that $f_{\lambda_2}(s\vert t)/f_{\lambda_1}(s\vert t)$ is increasing (decreasing) in $s$. Take $k_1(s)=I_{(-\infty,m_{\lambda_2}^{(1)}(t))}(s)$ and $k_2(s)=f_{\lambda_2}(s\vert t)/f_{\lambda_1}(s\vert t)$, where $I_A(\cdot)$ denotes the indicator function of set $A  \subseteq \Re$. Here $k_1(s)$ is decreasing in $s$ and $k_2(s)$ is increasing (decreasing) in $s$. Using Proposition 2.1, we get
	\begin{align*}
		\frac{1}{2}=&\int_{-\infty}^{\infty} k_1(s) k_2(s) f_{\lambda_1}(s\vert t)ds \leq (\geq) \left( \int_{-\infty}^{\infty} k_1(s) f_{\lambda_1}(s\vert t)ds \right) \left( \int_{-\infty}^{\infty} k_2(s) f_{\lambda_1}(s\vert t)ds\right)\\
		\implies &\int_{-\infty}^{m_{\lambda_2}^{(1)}(t)} f_{\lambda_2}(s\vert t)ds=\int_{-\infty}^{m_{\lambda_1}^{(1)}(t)} f_{\lambda_1}(s\vert t)ds=\frac{1}{2} \leq (\geq) \int_{-\infty}^{m_{\lambda_2}^{(1)}(t)} f_{\lambda_1}(s\vert t)ds\\
		\implies & m_{\lambda_1}^{(1)}(t)\leq (\geq)\, m_{\lambda_2}^{(1)}(t),
	\end{align*}
	establishing the assertion.
\end{proof}	 

\noindent
Under the assumptions of Lemma 3.1.1, for any fixed $t$, one may take
\begin{align}\label{eq:3.3}
	l^{(1)}(t)& =\inf_{\lambda\geq 0} m_{\lambda}^{(1)}(t) =m_0^{(1)}(t)\; (=\lim_{\lambda\to \infty} m_{\lambda}^{(1)}(t))\\
	\;\text{and}\;\;
	u^{(1)}(t)&=\sup_{\lambda\geq 0} m_{\lambda}^{(1)}(t)=\lim_{\lambda\to \infty} m_{\lambda}^{(1)}(t)\;(=m_0^{(1)}(t)),
\end{align}
while applying Theorem 3.1.1 and Corollary 3.1.1. \vspace*{2mm}

%			\\~\\\textbf{Corollary 3.1.3.} Let $\xi^{*}(t)\!=\!\max\{l^{(1)}(t),\min\{c_{0,1},u^{(1)}(t)\}\}$, where $l^{(1)}(\cdot)$ and $u^{(1)}(\cdot)$ are given by \eqref{eq:3.3}. Then
%		the estimator $\delta_{\xi^{*}}(\bold{X})\!=\!X_1-\xi^{*}(D)$ is nearest to $\theta_1$ than the BLEE $\delta_{1,c_{0,1}}(\bold{X})=X_1-c_{0,1}$ under the GPN criterion for all $\boldsymbol{\theta}\in \Theta_0$.\vspace*{2mm}

While applying Theorem 3.1.1 and Corollaries 3.1.1-3.12, for any fixed $t$, the commonly used choice for $(l^{(1)}(t),u^{(1)}(t))$ is given by $l^{(1)}(t)=\inf_{\lambda\geq 0} m_{\lambda}^{(1)}(t)$ and $u^{(1)}(t)=\sup_{\lambda\geq 0} m_{\lambda}^{(1)}(t)$. Now we will provide some applications of Theorem 3.1.1 and Corollaries 3.1.1-3.1.2.
\\~\\ \textbf{Example 3.1.1.} Let $\bold{X}=(X_1,X_2)$ follow a bivariate normal distribution with joint p.d.f. (3.1),
where, for known positive real numbers $\sigma_1$ and $\sigma_2$ and known $\rho \in (-1,1),$ the joint p.d.f. of $\bold{Z}=(Z_1,Z_2)=(X_1-\theta_1,X_2-\theta_2)$ is
$$ f(z_1,z_2) =\frac{1}{2 \pi \sigma_1 \sigma_2 \sqrt{1-\rho^2}} e^{-\frac{1}{2(1-\rho^2)}\left[\frac{z_1^2}{\sigma_1^2}-2 \rho \, \frac{z_1 z_2}{\sigma_1 \sigma_2}+\frac{z_2^2}{\sigma_2^2}\right]},\; \; \; \bold{z}=(z_1,z_2)\in \Re^2.$$ 
Consider estimation of $\theta_1$ under the GPN criterion with a general loss function (i.e., $L_1(\boldsymbol{\theta},a)=W(a-\theta_1),\;\boldsymbol{\theta}\in\Theta_0,\; a\in\Re$, where $W(0)=0$, $W(t)$ is strictly decreasing on $(-\infty,0)$ and strictly increasing on $(0,\infty)$). In this case, the PNLEE is $\delta_{1,PNLEE}(\bold{X})=X_1$. Also, for any fixed $t\in \Re$ and $\lambda\geq 0$, the conditional distribution of $Z_1$ given $D=t$ is $N\left((1-\alpha)(\lambda-t),\frac{(1-\rho^2)\sigma_1^2 \sigma_2^2}{\tau^2}\right), \text{ where } \tau^2=\sigma_1^2+\sigma_2^2-2\rho \sigma_1 \sigma_2$ and $\alpha=\frac{\sigma_2(\sigma_2-\rho\sigma_1)}{\tau^2}$.  \vspace*{2mm}

Here the restricted MLE is $\delta_{1,RMLE}(\bold{X})=X_1-(1-\alpha)\max\{0,-D\}$. Hwang and Peddada (\citeyear{MR1272076}) and Tan and Peddada (\citeyear{tan2000}) proposed alternative estimators for $\theta_1$ as $\delta_{1,HP}(\bold{X})=\min\{X_1,\alpha X_1+(1-\alpha)X_2\}=X_1-\max\{0,(\alpha-1)D\}$ and $\delta_{1,PDT}(\bold{X})=\min\{X_1,\beta(\alpha) X_1+(1-\beta(\alpha))X_2\}=X_1-(1-\beta(\alpha))\max\big\{0,-D\big\}$, receptively, where $\beta(\alpha)=\min\{1,\max\{0,\alpha\}\},\;\alpha\in \Re$.

\vspace*{2mm}

The median of the conditional conditional distribution of random variable $Z_1$, given $D=t$, is $m_{\lambda}^{(1)}(t)=(1-\alpha)(\lambda-t),\;t\in\Re,\;\lambda\geq 0.$ Clearly, $m_{\lambda}^{(1)}(t)$ is increasing in $\lambda\in [0,\infty)$, if $\alpha<1$ and decreasing in $\lambda\in [0,\infty)$, if $\alpha>1$. Thus, for $t\in \Re$, as in (3.4) and (3.5), we may take
\begin{align*}
	l^{(1)}(t)&=\inf_{\lambda\geq 0} m_{\lambda}^{(1)}(t)=\begin{cases}
		-(1-\alpha)t,\;\; &\alpha\leq 1\\
		-\infty,\;\; &\alpha>1
	\end{cases}\\
	\text{and}\quad
	u^{(1)}(t)&=\sup_{\lambda\geq 0} m_{\lambda}^{(1)}(t)=\begin{cases}
		\infty,\;\; &\alpha<1\\
		-(1-\alpha)t,\;\; &\alpha\geq 1
	\end{cases}.
\end{align*}
Consider the following cases:
\vspace*{2mm}

\noindent
\textbf{Case-I:}  $0\leq \alpha <1$
\vspace{2mm} 

In this case the Hwang and Peddada (\citeyear{MR1272076}) estimator $\delta_{1,HP}(\bold{X})$, Tan and Peddada (\citeyear{tan2000}) estimator $\delta_{1,PDT}(\bold{X})$ and the restricted MLE $\delta_{1,RMLE}(\bold{X})$ are the same and no improvement is possible over these estimators using our results.

We have $l^{(1)}(t)=-(1-\alpha)t
\;\text{and}\;
u^{(1)}(t)=\infty,\;t\in \Re.$ Use of Theorem 3.1.1 and Corollary 3.1.1 leads us to following conclusions:
\vspace*{2mm}

\noindent
\textbf{(i)}  The estimator $\delta_{1,RMLE}(\bold{X})=X_1-(1-\alpha) \max\{0,-D\}$ ($=\delta_{1,HP}(\bold{X})=\delta_{1,PDT}(\bold{X})$) is nearer to $\theta_1$ than the estimator $\delta_{1,PNLEE}(\bold{X})=X_1$.
\vspace*{2mm}

\noindent
\textbf{(ii)}  Let $\psi_{1,0}(\cdot)$ be such that $0< \psi_{1,0}(t)\leq -(1-\alpha)t,\; \forall \; t<0$, and $\psi_{1,0}(t)=0,\; \forall \; t\geq 0$. Then the estimator $\delta_{\psi_{1,0}}(\bold{X})=X_1-\psi_{1,0}(D)$ is nearer to $\theta_1$ than $\delta_{1,PNLEE}(\bold{X})$. In particular, for the choice,
$$\psi_{1,\nu}(t)=\begin{cases} -(1-\nu)t,\;\; & t\leq 0\\ 0,\;\; & t> 0 \end{cases},\quad \alpha\leq \nu <1,$$
the estimator $\delta_{\psi_{1,\nu}}(\bold{X})=X_1-\psi_{1,\nu}(D)$ is nearer to $\theta_1$ than $\delta_{1,PNLEE}(\bold{X})$.

\vspace*{1.5mm}

\noindent \textbf{Case-II:}  $\alpha_1=1$
\vspace{2mm} 

In this case $\delta_{1,RMLE}(\bold{X})=\delta_{1,HP}(\bold{X})=\delta_{1,PDT}(\bold{X})=\delta_{1,PNLEE}(\bold{X})=X_1$. Also $l^{(1)}(t)=u^{(1)}(t)=0,\; \forall\; t\in \Re$. Each of the above estimators are nearer to $\theta_1$ than any other location equivariant estimator.

\vspace*{1.5mm}

\noindent\textbf{Case-III}  $\alpha>1$
\vspace{2mm} 

In this case, $\delta_{1,RMLE}(\bold{X})=X_1-(1-\alpha)\max\{0,-D\}$, $\delta_{1,HP}(\bold{X})=X_1-(\alpha-1)\max\{0,D\}$ and $\delta_{1,PDT}(\bold{X})=\delta_{1,PNLEE}(\bold{X})=X_1$. Also, $l^{(1)}(t)=-\infty
\;\text{and}\;
u^{(1)}(t)=-(1-\alpha)t, \; \forall \; t\in \Re.$ We have the following consequences of Theorem 3.1.1 and Corollary 3.1.1:

\vspace*{2mm}

\noindent
\textbf{(i)}  Estimators $\delta_{1,RMLE}(\bold{X})$ and $\delta_{1,HP}(\bold{X})$ are nearer to $\theta_1$ than the estimator $\delta_{1,PDT}(\bold{X})\;(=\delta_{1,PNLEE}(\bold{X}))$.
\vspace*{2mm}

\noindent
\textbf{(ii)}  The Estimator $\delta_{1,HP}^{*}(\bold{X})=\alpha X_1 +(1-\alpha)X_2$ is nearer to $\theta_1$ than the estimator $\delta_{1,HP}(\bold{X})$.
\vspace*{2mm}

\noindent
\textbf{(iii)}  Let $\psi_{1,0}(\cdot)$ be such that $-(1-\alpha)t\leq \psi_{1,0}(t)< 0,\; \forall \; t<0$, and $\psi_{1,0}(t)=0,\; \forall \; t\geq 0$. Then the estimator $\delta_{\psi_{1,0}}(\bold{X})=X_1-\psi_{1,0}(D)$ is nearer to $\theta_1$ than $\delta_{1,PNLEE}(\bold{X})$. In particular, for the choice,
$$\psi_{1,\nu}(t)=\begin{cases} -(1-\nu)t,\;\; & t\leq 0\\ 0,\;\; & t\geq 0 \end{cases},\quad 1< \nu \leq \alpha,$$
the estimator $\delta_{\psi_{1,\nu}}(\bold{X})=X_1-\psi_{1,\nu}(D)$ is nearer to $\theta_1$ than $\delta_{1,PNLEE}(\bold{X})$.

\vspace*{2mm}

\noindent
\textbf{(iv)}  Let $\psi_{1,0}(\cdot)$ be such that $-(1-\alpha)t\leq \psi_{1,0}(t)\leq 0,\; \forall \; t<0$, and $\psi_{1,0}(t)=-(1-\alpha)t,\; \forall \; t\geq 0$. Then the estimator $\delta_{\psi_{1,0}}(\bold{X})=X_1-\psi_{1,0}(D)$ is nearer to $\theta_1$ than $\delta_{1,HP}(\bold{X})$. In particular, for the choice, 
$$\psi_{1,\nu}(t)=\begin{cases} -(1-\nu)t,\;\; & t\leq 0\\ -(1-\alpha)t,\;\; & t> 0 \end{cases},\quad 1< \nu \leq \alpha,$$
the estimator $$\delta_{\psi_{1,\nu}}(\bold{X})=\begin{cases} \nu X_1 +(1-\nu) X_2,\;\; & X_1\geq X_2\\ \alpha X_1 +(1-\alpha) X_2,\;\; & X_1<X_2 \end{cases}$$ is nearer to $\theta_1$ than $\delta_{1,HP}(\bold{X})$.

\vspace*{2mm}

\noindent\textbf{Case-IV}  $\alpha<0$
\vspace{2mm} 

Here $\delta_{1,RMLE}(\bold{X})=\delta_{1,HP}(\bold{X})=X_1-(1-\alpha)\max\{0,-D\}$ and $\delta_{1,PDT}(\bold{X})=X_1-\max\{0,-D\}$. Also, we have $l^{(1)}(t)=-(1-\alpha)t
\;\text{and}\;u^{(1)}(t)=\infty, \; \forall \; t\in \Re.$ The following observations are evident from Theorem 3.1.1 and Corollary 3.1.1:

\vspace*{2mm}

\noindent
\textbf{(i)}  Estimators $\delta_{1,RMLE}(\bold{X})\;(=\delta_{1,HP}(\bold{X}))$ and $\delta_{1,PDT}(\bold{X})$ are nearer to $\theta_1$ than the estimator $\delta_{1,PNLEE}(\bold{X})$.

%\vspace*{2mm}

\noindent
\textbf{(ii)} The estimator $\delta_{1,RMLE}(\bold{X})\;(=\delta_{1,HP}(\bold{X}))$ is nearer to $\theta_1$ than $\delta_{1,PDT}(\bold{X})$.
\vspace*{2mm}

\noindent
\textbf{(iii)}  Let $\psi_{1,0}(\cdot)$ be such that $0< \psi_{1,0}(t)\leq -(1-\alpha)t,\; \forall \; t<0$, and $\psi_{1,0}(t)=0,\; \forall \; t\geq 0$. Then the estimator $\delta_{\psi_{1,0}}(\bold{X})=X_1-\psi_{1,0}(D)$ is nearer to $\theta_1$ than $\delta_{1,PNLEE}(\bold{X})=X_1$. In particular, for the choice,
$$\psi_{1,\nu}(t)=\begin{cases} -(1-\nu)t,\;\; & t\leq 0\\ 0,\;\; & t\geq 0 \end{cases},\quad \alpha\leq \nu <0,$$
the estimator $$\delta_{\psi_{1,\nu}}(\bold{X})=\begin{cases} \nu X_1 +(1-\nu) X_2,\;\; &X_1> X_2 \\ X_1,\;\; & X_1 \leq X_2 \end{cases}$$ 
is nearer to $\theta_1$ than $\delta_{1,PNLEE}(\bold{X})=X_1$.

\vspace*{2mm}

\noindent
\textbf{(iv)}  Let $\psi_{1,0}(\cdot)$ be such that $-t\leq \psi_{1,0}(t)\leq -(1-\alpha)t,\; \forall \; t<0$, and $\psi_{1,0}(t)=0,\; \forall \; t\geq 0$. Then the estimator $\delta_{\psi_{1,0}}(\bold{X})=X_1-\psi_{1,0}(D)$ is nearer to $\theta_1$ than $\delta_{1,PDT}(\bold{X})$. In particular, for the choice, 
$$\psi_{1,\nu}(t)=\begin{cases} -(1-\nu)t,\;\; & t\leq 0\\ 0,\;\; & t\geq 0 \end{cases},\quad \alpha\leq  \nu < 0,$$
the estimator $$\delta_{\psi_{1,\nu}}(\bold{X})=\begin{cases} \nu X_1 +(1-\nu) X_2,\;\; & X_1> X_2\\ X_1,\;\; & X_1\leq X_2 \end{cases}$$ is nearer to $\theta_1$ than $\delta_{1,PDT}(\bold{X})$.

\vspace*{3mm}

\noindent
For the bivariate normal model, some of the above results have also been reported in Chang et al. (\citeyear{chang2017}, \citeyear{chang2020}) under specific $W(t)=\vert t\vert,\;t\in\Re$. The findings reported in the above example hold for any general $W(\cdot)$ such that $W(0)=0$, $W(t)$ is increasing in $(0,\infty)$ and decreasing in $(-\infty,0)$.
%Now, the GPN of the restricted MLE $\delta_{R}(\bold{X})=X_1-\psi_{R}(D)$, where $\psi_{R}(t)=\min\{0,\frac{-\sigma_1(\sigma_1-\rho\sigma_2)t}{\tau^2}\},\;t\in\Re$, relative to $\delta_{HP}(\bold{X})=X_1-\psi_{HP}(D)$, where $\psi_{HP}(t)=\max\{0,\frac{-\sigma_1(\sigma_1-\rho\sigma_2)t}{\tau^2}\},\;t\in\Re$, can be written as (clearly, $\psi_{R}(t)<\psi_{HP}(t),\;t\in \Re$ )
%\begin{align*}
%	GPN(\delta_{R},\delta_{HP};\theta_1)&=P_{\boldsymbol{\theta}}[W(Z_1-\psi_{R}(D)) <W(Z_1-\psi_{HP}(D))],\;\; \boldsymbol{\theta}\in\Theta_0,\\
%	&=\int_{-\infty}^{\infty} P_{\boldsymbol{\theta}}[W(Z_1-\psi_{R}(t)) <W(Z_1-\psi_{HP}(t))\vert D=t] f_D(t) dt\\
%	&\geq \int_{-\infty}^{\infty} P_{\boldsymbol{\theta}}[Z_1\leq \psi_{R}(t)\vert D=t] f_D(t) dt\\
%		&= \int_{-\infty}^{\infty} P_{\boldsymbol{\theta}}[Z_1\leq \min\bigg\{0,\frac{-\sigma_1(\sigma_1-\rho\sigma_2)t}{\tau^2}\bigg\}\vert D=t] f_D(t) dt
%\end{align*}
%Since, for fix $t$, $\psi_{R}(t)<\psi_{HP}(t)$, then $z_1\leq \psi_{R}(t)$ implies that $W(z_1-\psi_{R}(t)) <W(z_1-\psi_{HP}(t))$.

%\vspace*{3mm}
\vspace*{3mm}

\noindent
\textbf{Example 3.1.2.} Let $X_1$ and $X_2$ be independent random variables with a joint p.d.f. $f(x_1-\theta_1,x_2-\theta_2)= \frac{1}{\sigma_1 \sigma_2}\, e^{-\frac{x_1-\theta_1}{\sigma_1}}e^{-\frac{x_2-\theta_2}{\sigma_2}},\;x_1>\theta_1,\;x_2>\theta_2,\;(\theta_1,\theta_2)\in\Theta_0$, where $\sigma_1$ and $\sigma_2$ are known positive constants. 
%Here, for every fixed $t\in \Re$ and $\lambda\geq 0$, $f(s,s+t-\lambda)/f(s,s+t)$ is increasing in $s\in[0,\infty)$. \vspace*{1.5mm}
Here the PNLEE is $\delta_{1,PNLEE}(\bold{X})=X_1-\sigma_1 \ln(2)$ and the restricted MLE is $\delta_{1,RMLE}(\bold{X})=\min\{X_1,X_2\}=X_1-\max\{0,-D\}$.

\vspace*{2mm}

For any fixed $t\in \Re$, the conditional p.d.f. of $Z_1$ given $D=t$ is 
\begin{align*}
	f_{\lambda}(s\vert t)&=\frac{\sigma_1+\sigma_2}{\sigma_1 \sigma_2}\, e^{-\left(\frac{1}{\sigma_1}+\frac{1}{\sigma_2}\right)(s-\max\{-t+\lambda,0\})},  \; s\geq \max\{-t+\lambda,0\}
	\\ \text{and} \;\;\;\; m_{\lambda}^{(1)}(t)&=\max\{-t+\lambda,0\}+\frac{\sigma_1\sigma_2}{\sigma_1+\sigma_2} \ln(2),\;\; t\in \Re.
\end{align*}
Clearly, for every fixed $t\in \Re$, the median $m_{\lambda}^{(1)}(t)$ is an increasing function of $\lambda\in[0,\infty)$ (this also follows from Lemma 3.1.1 as, for every fixed $\lambda\geq 0$ and $t\in \Re$, $f(s,s+t-\lambda)/f(s,s+t)$ is increasing in $s\in (0,\infty)$). Thus, we may take
$l^{(1)}(t)=\!\inf_{\lambda\geq 0} m_{\lambda}^{(1)}(t)= \max\{0,-t\}+\frac{\sigma_1 \sigma_2}{\sigma_1+\sigma_2}\ln(2)
\text{ and }
u^{(1)}(t)\! 
=\sup_{\lambda\geq 0} m_{\lambda}^{(1)}(t)=\infty,\;t\in\Re.$ 
\vspace*{2mm}

\noindent The following conclusions immediately follow from Theorem 3.1.1 and Corollary 3.1.1:
\\~\\\textbf{(i)}  The estimator $\delta_{1,PNLEE}^{*}(\bold{X})=\min\big\{X_2-\frac{\sigma_1\sigma_2}{\sigma_1+\sigma_2}\ln(2),X_1-\sigma_1\ln(2)\big\}$ is nearer to $\theta_1$ than the PNLEE $\delta_{1,PNLEE}(\bold{X})$.
\\~\\ \textbf{(ii)}  The estimator
$\delta_{1,RMLE}^{*}(\bold{X})=\min\{X_1,X_2\}-\frac{\sigma_1 \sigma_2}{\sigma_1+\sigma_2}\ln(2)$ is nearer to $\theta_1$ than the restricted MLE $\delta_{1,RMLE}$.
\\~\\ \textbf{(iii)}  Let $\psi_{1,0}(t)$ be such that $\sigma_1 \ln(2) < \psi_{1,0}(t) \leq -t + \frac{\sigma_1 \sigma_2}{\sigma_1+\sigma_2}\ln(2),\; \forall\; t\leq \frac{-\sigma_1^2}{\sigma_1+\sigma_2}\ln(2)$, and $\psi_{1,0}(t)=\sigma_1 \ln(2),\; \forall\; t>\frac{-\sigma_1^2}{\sigma_1+\sigma_2}\ln(2)$. Then the estimator $\delta_{\psi_{1,0}}(\bold{X})=X_1-\psi_{1,0}(D)$ is nearer to $\theta_1$ than the PNLEE $\delta_{1,PNLEE}(\bold{X})=X_1-\sigma_1 \ln(2)$.
\\~\\ \textbf{(iv)}  Let $\psi_{1,0}(t)$ be such that $-t \leq \psi_{1,0}(t) \leq -t + \frac{\sigma_1 \sigma_2}{\sigma_1+\sigma_2}\ln(2),\; \forall\; t\leq 0$, and $0 \leq \psi_{1,0}(t) \leq \frac{\sigma_1 \sigma_2}{\sigma_1+\sigma_2}\ln(2),\; \forall\; t> 0$. Then the estimator $\delta_{\psi_{1,0}}(\bold{X})=X_1-\psi_{1,0}(D)$ is nearer to $\theta_1$ than the restricted MLE $\delta_{1,RMLE}(\bold{X})=\min\{X_1,X_2\}$.

%The restricted MLE of $\theta_1$ is $\delta_{1,R}(\bold{X})\!=\!\min\{X_1,X_2\}\!=\!X_1-\max\{0,-D\}$. Using Theorem 3.1.1, it follow that the estimator
%$\delta_{1,R}^{*}(\bold{X})=\min\{X_1,X_2\}-\frac{\sigma_1 \sigma_2}{\sigma_1+\sigma_2}\ln(2)$ is nearest to $\theta_1$ than the restricted MLE $\delta_{1,R}$ under the GPN criterion.
%\vspace*{1mm}

%Under the GPN criterion, the unrestricted PNEE of $\theta_1$ is $\delta_{1,m_{0,1}}(\bold{X})=X_1-\sigma_1\ln(2)$ (i.e., $m_{0,1}=\sigma_1\ln(2)$).
%Using Corollary 3.1.2, it follow that the the estimator
%$\delta_{1,m_{0,1}}^{*}(\bold{X})=\min\big\{X_2-\frac{\sigma_1\sigma_2}{\sigma_1+\sigma_2}\ln(2),X_1-\sigma_1\ln(2)\big\}$ is closet to $\theta_1$ than the PNEE $\delta_{1,m_{0,1}}(\bold{X})=X_1-\sigma_1\ln(2)$ under the GPN criterion.
%\vspace*{1.5mm}

%Also, under the mean squared error (MSE) criterion, the unrestricted BLEE of $\theta_1$ is $\delta_{1,c_{0,1}}(\bold{X})=X_1-\sigma_1$ (i.e., $c_{0,1}=\sigma_1$).
%Using Corollary 3.1.3, it follow that the the estimator
%$\delta_{1,c_{0,1}}^{*}(\bold{X})=\min\big\{X_2-\frac{\sigma_1\sigma_2}{\sigma_1+\sigma_2}\ln(2),X_1-\sigma_1\big\}$ is nearest to $\theta_1$ than the BLEE $\delta_{1,c_{0,1}}(\bold{X})=X_1-\sigma_1$ under the GPN criterion.
\vspace*{2.5mm}

\noindent
%	\textbf{Remark 3.1.1.} For the exponential model considered in Example 3.1.2, Vijayasree et al. (\citeyear{MR1345425}) used the minimizing mean squared error (MSE) criterion and obtain the dominating estimator over the BLEE $\delta_{1,c_{0,1}}(\bold{X})=X_1-\sigma_1$ which is $ \delta_{MSE}(\bold{X})=\min\big\{X_2-\frac{\sigma_1\sigma_2}{\sigma_1+\sigma_2},X_1-\sigma_1\big\}.$
It is worth mentioning that findings of Theorem 3.1.1 and Corollaries 3.1.1-3.1.2, and hence those of Examples 3.1.1 and 3.1.2, hold under any general loss function $L_1(\boldsymbol{\theta},a)=W(a-\theta_1),\;\boldsymbol{\theta}\in\Theta_0,\; a\in\mathcal{A}=\Re,$     
where $W:\Re\rightarrow [0,\infty)$ is such that $W(0)=0$, $W(t)$ is strictly decreasing on $(-\infty,0)$ and strictly increasing on $(0,\infty)$.

\subsection{\textbf{Estimation of the Larger Location Parameter $\theta_2$}} \label{3.2}

Consider estimation of the larger location parameter $\theta_2$ under the GPN criterion with the general loss function $L_2(\boldsymbol{\theta},a)=W(a-\theta_2),\;\boldsymbol{\theta}\in\Theta_0,\; a\in\mathcal{A}= \Re$, when it is known apriori that $\boldsymbol{\theta}\in\Theta_0$. The form of any location equivariant estimator of $\theta_2$ is $\delta_{\psi}(\bold{X})=X_2-\psi(D),$ for some function $\psi:\,\Re\rightarrow \Re$, where $D=X_2-X_1$.\vspace*{2mm}

Let $Z_2=X_2-\theta_2$, $\lambda=\theta_2-\theta_1$ and $f_D(t\vert \lambda)$ be the p.d.f. of r.v. $D=X_2-X_1$. Let $\delta_{\xi}(\bold{X})=X_2-\xi(D)$ and $\delta_{\psi}(\bold{X})=X_2-\psi(D)$ be two location equivariant estimators of $\theta_2$. Then, the GPN of $\delta_{\xi}(\bold{X})$ relative to $\delta_{\psi}(\bold{X})$ is given by
\begin{align*}
	GPN(\delta_{\xi},\delta_{\psi};\boldsymbol{\theta})
	=\int_{-\infty}^{\infty} g_{2,\lambda}(\xi(t),\psi(t),t) f_D(t\vert \lambda) dt,\;\; \lambda\geq 0,
\end{align*}
where, for $\lambda \geq 0$, 

\begin{align*}g_{2,\lambda}(\xi(t),\psi(t),t)&= P_{\boldsymbol{\theta}}[W(Z_2-\xi(t)) <W(Z_2-\psi(t))\vert D=t]
	\\&\qquad +\frac{1}{2}P_{\boldsymbol{\theta}}[W(Z_2-\xi(t))\! =\!W(Z_2-\psi(t))\vert D=t].\end{align*}
Let $m_{\lambda}^{(2)}(t)$ denote the median of the conditional distribution of $Z_2$ given $D=t$, where $\lambda \geq 0$ and $t$ belongs to the support of r.v. $D$. For any fixed $t\in \Re$, the conditional p.d.f. of $Z_2$ given $D=t$ is $f_{\lambda}(s\vert t)=\frac{f(s-t+\lambda,s)}{f_D(t\vert \lambda)}$ and $f_D(t\vert \lambda)=\int_{-\infty}^{\infty}f(y-t
+\lambda,y)dy$, $\lambda\geq 0$. Thus $\int_{-\infty}^{m_{\lambda}^{(2)}(t)} \, f(s-t+\lambda,s)ds=\frac{1}{2} \int_{-\infty}^{\infty}f(s-t+\lambda,s)ds$. For any fixed $t$ and $\lambda\geq 0$, using Lemma 3.1, we have $g_{2,\lambda}(\xi(t),\psi(t),t)> \frac{1}{2}$, provided $\psi(t)<\xi(t)\leq m_{\lambda}^{(2)}(t)$ or if $ m_{\lambda}^{(2)}(t)\leq \xi(t)<\psi(t)$. Also, for any fixed $t$, $g_{2,\lambda}(\psi(t),\psi(t),t)= \frac{1}{2},\;\forall \;\lambda\geq 0$. Now using arguments similar to the ones used in proving Theorem 3.1.1, we get the following results.  \vspace*{2mm} 

\noindent
\textbf{Theorem 3.2.1.} Let $\delta_{\psi}(\bold{X})=X_2-\psi(D)$ be a location equivariant estimator of $\theta_2$. Let $l^{(2)}(t)$ and $u^{(2)}(t)$ be functions such that $l^{(2)}(t)\leq m_{\lambda}^{(2)}(t)\leq u^{(2)}(t), \;\forall\;\lambda\geq 0$ and any $t$. For any fixed $t$, define $\psi^{*}(t)\!=\!\max\{l^{(2)}(t),\min\{\psi(t),u^{(2)}(t)\}\}$. Let $\delta_{\psi^{*}}(\bold{X})\!=\!X_2-\psi^{*}(D)$. Then, 
$GPN(\delta_{\psi^{*}},\delta_{\psi};\boldsymbol{\theta})> \frac{1}{2},\; \forall \; \boldsymbol{\theta}\in \Theta_0,$ provided $P_{\boldsymbol{\theta}}[l^{(2)}(D)\leq \psi(D)\leq u^{(2)}(D)]<1,\;\forall\; \boldsymbol{\theta}\in \Theta_0$.
\vspace*{2mm}

\noindent
\textbf{Corollary 3.2.1.} Let $\delta_{\psi}(\bold{X})=X_2-\psi(D)$ be a location equivariant estimator of $\theta_2$ such that $P_{\boldsymbol{\theta}}[l^{(2)}(D)\leq \psi(D)\leq u^{(2)}(D)]<1,\;\forall\; \boldsymbol{\theta}\in \Theta_0$, where $l^{(2)}(\cdot)$ and $u^{(2)}(\cdot)$ are as defined in Theorem 3.2.1. Let $\psi_{2,0}:\Re\rightarrow \Re$ be such that $\psi(t)< \psi_{2,0}(t) \leq l^{(2)}(t)$, whenever $\psi(t)< l^{(2)}(t)$, or $u^{(2)}(t)\leq \psi_{2,0}(t)< \psi(t)$, whenever $u^{(2)}(t)< \psi(t)$. Also let $\psi_{2,0}(t)=\psi(t)$, whenever $l^{(2)}(t)\leq \psi(t) \leq u^{(2)}(t)$. Let $\delta_{\psi_{2,0}}(\bold{X})=X_2-\psi_{2,0}(D)$. Then, 
$GPN(\delta_{\psi_{2,0}},\delta_{\psi};\boldsymbol{\theta})> \frac{1}{2},\; \forall \; \boldsymbol{\theta}\in \Theta_0$.\vspace*{2mm}

The following corollary provides improvements over the PNLEE $\delta_{2,PNLEE}(\bold{X})=X_2-m_{0,2}$, under the restricted parameter space $\Theta_0$
\\~\\\textbf{Corollary 3.2.2.} Let $\xi^{*}(t)\!=\!\max\{l^{(2)}(t),\min\{m_{0,2},u^{(2)}(t)\}\},\; t\in \Re$, where $l^{(2)}(\cdot)$ and $u^{(2)}(\cdot)$ are as defined in Theorem 3.2.1. Let $\delta_{\xi^{*}}(\bold{X})\!=\!X_2-\xi^{*}(D)$. Then,
$GPN(\delta_{\xi^{*}},\delta_{2,PNLEE};\boldsymbol{\theta})> \frac{1}{2},\; \forall \; \boldsymbol{\theta}\in \Theta_0,$ provided $P_{\boldsymbol{\theta}}[l^{(2)}(D)\leq m_{0,2}\leq u^{(2)}(D)]<1,\;\forall\; \boldsymbol{\theta}\in \Theta_0$.
\vspace*{2mm}

The following lemma describes the behaviour of $m_{\lambda}^{(2)}(t)$, for any fixed $t$. The proof of the lemma, being similar to the proof of Lemma 3.1.1, is skipped.
\vspace*{2mm}		

\noindent \textbf{Lemma 3.2.1.} If, for every fixed $\lambda\geq 0$ and $t$, $f(s-t+\lambda,s)/f(s-t,s)$ is increasing (decreasing) in $s$ (wherever the ratio is not of the form $0/0$), then, for every fixed $t$, $m_{\lambda}^{(2)}(t)$ is an increasing (decreasing) function of $\lambda\in [0,\infty)$.
%\vspace*{2mm}

\noindent
Under the assumptions of Lemma 3.2.1, one may take, for any fixed $t$, \begin{align}\label{eq:3.3}
	l^{(2)}(t)& =\inf_{\lambda\geq 0} m_{\lambda}^{(2)}(t) =m_0^{(2)}(t)\; (=\lim_{\lambda\to \infty} m_{\lambda}^{(2)}(t))\\
	\;\text{and}\;\;
	u^{(2)}(t)&=\sup_{\lambda\geq 0} m_{\lambda}^{(2)}(t)=\lim_{\lambda\to \infty} m_{\lambda}^{(2)}(t)\;(=m_0^{(2)}(t)),
\end{align}
while applying Theorem 3.2.1 and Corollary 3.2.1. \vspace*{2mm}

As in Section 3.1, we will now apply Theorem 3.2.1 and Corollaries 3.2.1-3.2.2 to estimation of the larger location parameter $\theta_2$ in probability models considered in Examples 3.1.1-3.1.2.
\\~\\ \textbf{Example 3.2.1.} Let $\bold{X}=(X_1,X_2)$ have a bivariate normal distribution as described in Example 3.1.1. Consider estimation of $\theta_2$ under the GPN criterion with a general loss function $L_2(\boldsymbol{\theta},a)=W(a-\theta_2),\;\boldsymbol{\theta}\in\Theta_0,\; a\in\Re$, where $W(0)=0$, $W(t)$ is strictly decreasing on $(-\infty,0)$ and strictly increasing on $(0,\infty)$. Here, for any fixed $t\in \Re$, $Z_2 \vert D=t\sim N\left(\alpha(t-\lambda),\frac{(1-\rho^2)\sigma_1^2 \sigma_2^2}{\tau^2}\right), \text{ where } \tau^2=\sigma_1^2+\sigma_2^2-2\rho \sigma_1 \sigma_2$ and $\alpha=\frac{\sigma_2(\sigma_2-\rho\sigma_1)}{\tau^2}$. Thus, for $\lambda\geq 0$, $m_{\lambda}^{(2)}(t)=\alpha(t-\lambda),\;t\in \Re,$ and as in (3.6) and (3.7), we may take  
\begin{align*}
	l^{(2)}(t)=\begin{cases}
		\alpha t,\;\; &\alpha\leq 0\\
		-\infty,\;\; &\alpha>0
	\end{cases}\;\;
	\text{and}\quad
	u^{(2)}(t)=\begin{cases}
		\infty,\;\; &\alpha\leq 0\\
		\alpha t,\;\; &\alpha>0
	\end{cases}.
\end{align*}

The unrestricted PNLEE of $\theta_2$ is $\delta_{2,PNLEE}(\bold{X})=X_2$ (as $m_{0,2}=0$)
and the restricted MLE of $\theta_2$ is $\delta_{2,RMLE}(\bold{X})=X_2+\alpha \max\{0,-D\}$. Hwang and Peddada (\citeyear{MR1272076}) and Tan and Peddada (\citeyear{tan2000}) proposed alternative estimators for $\theta_2$ as $\delta_{2,HP}(\bold{X})=X_2+\max\{0,-\alpha D\}$ and $\delta_{2,PDT}(\bold{X})=X_2+\beta(\alpha)\max\big\{0,-D\big\}$, respectively, where $\beta(\alpha)=\min\{1,\max\{0,\alpha\}\},\;\alpha\in \Re$.
\vspace*{2mm}

Using Corollary 3.2.2, we conclude that, under the GPN criterion, the estimator 

$$\delta_{2,PNLEE}^{*}(\bold{X})=\delta_{2,RMLE}^{*}(\bold{X})=\begin{cases} \min\Big\{X_2, \alpha X_1+(1-\alpha) X_2\Big\},  & \alpha \leq 0
	\\~\\
	\max\Big\{X_2, \alpha X_1+(1-\alpha) X_2\Big\},  & \alpha>0 \end{cases}$$ 
is nearer to $\theta_2$ than $\delta_{2,PNLEE}(\bold{X})=X_2.$ In the same line as in Example 3.1.1, estimators dominating over $\delta_{2,HP}(\bold{X})$ and $\delta_{2,PDT}(\bold{X})$ can be obtained in certain situations.\vspace*{2mm}

\noindent 
\textbf{Example 3.2.2.} Let $X_1$ and $X_2$ be independent exponential random variables as considered in Example 3.1.2. Consider estimation of $\theta_2$ under the GPN criterion. Here the PNLEE of $\theta_2$ is $\delta_{2,PNLEE}(\bold{X})=X_2-\sigma_2\ln(2)$ and the restricted MLE of $\theta_2$ is $\delta_{2,RMLE}(\bold{X})=X_2$. Also, for any fixed $\lambda\geq 0$ and $t\in\Re$, the conditional p.d.f. of $Z_2$ given $D=t$ is
$$f_{\lambda}(s\vert t)=\frac{\sigma_1+\sigma_2}{\sigma_1 \sigma_2}\, e^{-\left(\frac{1}{\sigma_1}+\frac{1}{\sigma_2}\right)(s-\max\{t-\lambda,0\})}, \text{ if} \; s\geq \max\{t-\lambda,0\}.$$
Consequently,
$$ m_{\lambda}^{(2)}(t)= \max\{0,t-\lambda\}+\frac{\sigma_1 \sigma_2}{\sigma_1+\sigma_2}\ln(2),\; t\in \Re,\; \; \lambda \geq 0,$$
and, as in (3.6) and (3.7), we may take
$l^{(2)}(t)=\frac{\sigma_1 \sigma_2}{\sigma_1+\sigma_2}\ln(2),
\; t\in \Re,$ and $
u^{(2)}(t)= \max\{t,0\}+\frac{\sigma_1 \sigma_2}{\sigma_1+\sigma_2}\ln(2),\; t\in \Re.$
\vspace*{1.5mm}

The following conclusions are evident from Theorem 3.2.1 and Corollary 3.2.1:
\\~\\ \textbf{(i)} The estimator $\delta_{2,RMLE}^{*}(\bold{X})=X_2-\frac{\sigma_1 \sigma_2}{\sigma_1+\sigma_2} \ln(2)$ is nearer to $\theta_2$ than the restricted MLE $\delta_{2,RMLE}(\bold{X})=X_2$.
\\~\\ \textbf{(ii)} The estimator  \begin{align*}
	\delta_{2,PNLEE}^{*}(\bold{X})=\begin{cases} X_2-\frac{\sigma_1\sigma_2}{\sigma_1+\sigma_2}\ln(2),& \text{if } X_2< X_1\\ X_1-\frac{\sigma_1\sigma_2}{\sigma_1+\sigma_2}\ln(2),& \text{if }X_1\leq X_2< X_1+\frac{\sigma_2^2}{\sigma_1+\sigma_2}\ln(2)\\ X_2-\sigma_2\ln(2),&\text{if } X_2\geq X_1+\frac{\sigma_2^2}{\sigma_1+\sigma_2}\ln(2) \end{cases}
\end{align*} 
is nearer to $\theta_2$ than $\delta_{2,PNLEE}(\bold{X})$.
\\~\\ \textbf{(iii)}    Let $\psi_{2,0}(t)$ be such that $\max\{0,t\} +\frac{\sigma_1 \sigma_2}{\sigma_1 +\sigma_2} \ln(2) \leq \psi_{2,0}(t) < \sigma_2 \ln(2),\; \forall\; t\leq \frac{\sigma_2^2}{\sigma_1+\sigma_2}\ln(2)$, and $\psi_{2,0}(t)=\sigma_2 \ln(2),\; \forall\; t>\frac{\sigma_2^2}{\sigma_1+\sigma_2}\ln(2)$. Then the estimator $\delta_{\psi_{2,0}}(\bold{X})=X_2-\psi_{2,0}(D)$ is nearer to $\theta_2$ than the PNLEE $\delta_{2,PNLEE}(\bold{X})=X_2-\sigma_2 \ln(2)$.

%Also, under the mean squared error (MSE) criterion, the unrestricted BLEE of $\theta_2$ is $\delta_{2,c_{0,2}}(\bold{X})=X_2-\sigma_2$ (i.e., $c_{0,2}=\sigma_2$).
%Using Corollary 3.3.3, $GPN(\delta_{2,c_{0,2}}^{*},\delta_{2,c_{0,2}};\boldsymbol{\theta})\geq \frac{1}{2},\; \forall \; \boldsymbol{\theta}\in \Theta_0,$
%where
%\begin{align*}
%	\delta_{2,c_{0,2}}^{*}(\bold{X})=\begin{cases} X_2-\frac{\sigma_1\sigma_2}{\sigma_1+\sigma_2}\ln(2),& X_2< X_1\\ X_1-\frac{\sigma_1\sigma_2}{\sigma_1+\sigma_2}\ln(2),& X_1\leq X_2< X_1+\sigma_2-\frac{\sigma_2^2}{\sigma_1+\sigma_2}\ln(2)\\ X_2-\sigma_2,& X_2\geq X_1+\sigma_2-\frac{\sigma_2^2}{\sigma_1+\sigma_2}\ln(2) \end{cases}.
%\end{align*}

\section{\textbf{Improved Estimators for Restricted Scale Parameters}} \label{4}

Let $\bold{X}=(X_1,X_2)$ be a random vector having a joint p.d.f.
\begin{equation}\label{eq:3.1}
	f_{\boldsymbol{\theta}}(x_1,x_2)= \frac{1}{\theta_1 \theta_2}	f\left(\frac{x_1}{\theta_1},\frac{x_2}{\theta_2}\right),\; \; \;(x_1,x_2)\in \Re^2, 
\end{equation} 
where $f(\cdot,\cdot) $ is a specified Lebesgue p.d.f. and, for $\Re_{++}=(0,\infty)\times (0,\infty)$, $\boldsymbol{\theta}=(\theta_1,\theta_2)\in \Theta_0=\{(t_1,t_2)\in\Re_{++}^2:t_1 \leq t_2\}$ is the vector of unknown restricted scale parameters. For the sake of simplicity, throughout this section, we assume that the distributional support of $\bold{X}=(X_1,X_2)$ is a subset of $\Re_{++}^2$.  Generally, $\bold{X}=(X_1,X_2)$ would be a minimal-sufficient statistic based on a bivariate random sample or two independent random samples.  \vspace*{1.5mm}

For estimation of $\theta_i$, under the restricted parameter space $\Theta_0$, we consider the GPN criterion with the loss function $L_i(\boldsymbol{\theta},a)=W(\!\frac{a}{\;\theta_i}\!),\; \boldsymbol{\theta}\in\Theta_0,\; \; a\in \mathcal{A}=\Re_{++}=(0,\infty),\;i=1,2,$     
where $W:\Re_{++}\rightarrow [0,\infty)$ is a function such that $W(1)=0$, $W(t)$ is strictly decreasing on $(0,1)$ and strictly  increasing on $(1,\infty)$. Every time the word "general loss function" is used in this section, it refers to the loss function as defined above.
\vspace*{1.5mm}

The problem of estimating $\theta_i$, under the restricted parameter space $\Theta_0$ and under the GPN criterion with a general loss function defined above, is invariant under the group of transformations $\mathcal{G}\!=\!\{g_b:\,b\in(0,\infty)\},$ where $g_b(x_1,x_2)\!=\!(b\,x_1,b\,x_2),\; (x_1,x_2)\in\Re^2,\;b\in(0,\infty)$. Any scale equivariant estimator of $\theta_i$ has the form $$\delta_{\psi}(\bold{X})\!=\!\psi(D) X_i,$$
for some function $\psi:\,\Re_{++}\rightarrow \Re_{++}\,,\;i=1,2,$ where $D=\frac{X_2}{X_1}$. Let $f_D(t\vert \lambda)$ be the p.d.f. of r.v. $D=\frac{X_2}{X_1}$, where $\lambda=\frac{\theta_2}{\theta_1}\in[1,\infty)$. Note that the distribution of $D$ depends on $\boldsymbol{\theta}\in\Theta_0$ only through $\lambda=\frac{\theta_2}{\theta_1}\in[1,\infty)$. Exploiting the prior information of order restriction on parameters $\theta_1$ and $\theta_2$ ($\theta_1\leq \theta_2$), our aim is to obtain estimators that are nearer to $\theta_i,\,i=1,2$.\vspace*{2mm}

\noindent
The following lemma, whose proof is similar to that of Lemma 3.1, will play an important role in proving the main results of this section.
\\~\\	 \textbf{Lemma 4.1} Let $Y$ be a positive r.v. ($Pr(Y>0)=1$) having the Lebesgue p.d.f. and let $m_Y>0$ be the median of $Y$. Let $W:\Re_{++}\rightarrow [0,\infty)$ be a function such that $W(1)=0$, $W(t)$ is strictly decreasing on $(0,1)$ and strictly increasing on $(1,\infty)$. Then, for $0<c_1< c_2\leq m_{Y}$ or $0<m_Y\leq c_2<c_1$,
$GPN= P\left[W\left(\frac{Y}{c_2}\right)<W\left(\frac{Y}{c_1}\right)\right]+\frac{1}{2} P\left[W\left(\frac{Y}{c_2}\right)=W\left(\frac{Y}{c_1}\right)\right]> \frac{1}{2}$.

%	\noindent	\textbf{Lemma 4.1} Let $X$ be a r.v. with density $\frac{1}{\theta}f(\frac{x}{\theta}),\;x>0$, where $\theta>0$ is unknown scale parameter. Then,
%	$\delta_1=\frac{X}{m}$ is better than the estimator $\delta_2=\frac{X}{c}$ under the GPN criterion with the general loss function, where $m>0$ is the unique median of the r.v. $Z=\frac{X}{\theta}$ and $c>0$ is any constant.
%\begin{proof}.  For $Z=\frac{X}{\theta}$, GPN of $\delta_1$ relative to $\delta_2$ is
%	$$GPN(\delta_1,\delta_2;\theta)=P\left[W\left(\frac{Z}{m}\right) <W\left(\frac{Z}{c}\right)\right]+\frac{1}{2} P\left[W\left(\frac{Z}{m}\right) =W\left(\frac{Z}{c}\right)\right].$$	
%	Now, we have the following three cases
%		\\Case 1: If $m=c$, then $GPN(\delta_1,\delta_2;\theta)= \frac{1}{2}, \; \forall \; \theta\in \Theta$. \vspace*{1.2mm}

%		\noindent Case 2: If $m<c$, then $z\leq  m$ implies $W(\frac{z}{m})<W(\frac{z}{c})$ and we get,
%		$$GPN(\delta_1,\delta_2;\theta)=P\left[W\left(\frac{Z}{m}\right) <W\left(\frac{Z}{c}\right)\right]\geq P[Z\leq m]= \frac{1}{2}, \; \forall \; \theta\in \Theta.$$
%		\\Case 3: If $m>c$, then $z\geq  m$ implies $W(\frac{z}{m})<W(\frac{z}{c})$ and we get,
%		$$GPN(\delta_1,\delta_2;\theta)=P\left[W\left(\frac{Z}{m}\right) <W\left(\frac{Z}{c}\right)\right]\geq P[Z\geq  m]= \frac{1}{2}, \; \forall \; \theta\in \Theta.$$
%	\end{proof}
\vspace*{2mm}

Note that, in the unrestricted case (parameter space $\Theta\!=\!\Re_{++}^2$), the problem of estimating the scale parameter $\theta_i,\,i=1,2$, under the GPN criterion with a general loss function, is invariant under the multiplicative group of transformations $\mathcal{G}_0=\{g_{c_1,c_2}:\,(c_1,c_2)\in\Re_{++}^2\},$ where $g_{c_1,c_2}(x_1,x_2)=(c_1x_1,c_2x_2),\; (x_1,x_2)\in\Re^2,\;(c_1,c_2)\in\Re_{++}^2$. Any scale equivariant estimator is of the form $\delta_{i,c}(\bold{X})=cX_i,\;c\in \Re_{++},\; i=1,2.$ Using Lemma 4.1, the unrestricted PNSEE of $\theta_i$ is $\delta_{i,PNSEE}(\bold{X})=\frac{X_i}{m_{0,i}}$, where $m_{0,i}>0$ is the median of the r.v. $Z_i=\frac{X_i}{\theta_i},\,i=1,2$. \vspace*{2mm}

% Also, under the unrestricted parameter space (i.e., $\Theta=\Re_{++}^2$), the problem of estimating the scale parameter $\theta_i$, under the loss function $L_i(a,\boldsymbol{\theta})=W\left(\frac{a}{\theta_i}\right),\;\boldsymbol{\theta}\in\Theta_0, \;a\in\Re_{++}, \,i=1,2$, is invariant under the multiplicative group of transformations $\mathcal{G}_0$. The unrestricted best scale equivariant estimator (BSEE) of $\theta_i$ is $\delta_{i,BSEE}(\bold{X})=c_{0,i}X_i$, where $c_{0,i}$ is the unique solution of the equation $\int_{0}^{\infty} \int_{0}^{\infty}\, s_iW'(cs_i) \,f(s_1,s_2)\,ds_1\,ds_2=0,\;i=1,2.$\vspace*{2mm}

In the following subsections, we consider component-wise estimation of order restricted scale parameters $\theta_1$ and $\theta_2$, under the GPN criterion with a general loss function, and derive some general results. Applications of main results are illustrated through various examples dealing with specific probability models. \vspace{2mm}

\subsection{\textbf{Estimation of The Smaller Scale Parameter $\theta_1$}} \label{4.1}

Define $Z_1=\frac{X_1}{\theta_1}$ and $\lambda=\frac{\theta_2}{\theta_1}$, so that $\lambda\geq 1$. Let $f_D(t\vert \lambda)$ be the p.d.f. of r.v. $D=\frac{X_2}{X_1}$. Let $\delta_{\psi}(\bold{X})= \psi(D) X_1$ and $\delta_{\xi}(\bold{X})=\xi(D) X_1$ be two scale equivariant estimators of $\theta_1$, where $\psi:\,\Re_{++}\rightarrow \Re_{++}$ and $\xi:\,\Re_{++}\rightarrow \Re_{++}$ are specified functions. Then, the GPN of $\delta_{\xi}(\bold{X})=\xi(D)X_1$ relative to $\delta_{\psi}(\bold{X})=\psi(D)X_1$ is given by
\begin{align*}
	GPN(\delta_{\xi},\delta_{\psi};\boldsymbol{\theta})=\int_{0}^{\infty} g_{1,\lambda}(\xi(t),\psi(t),t) f_D(t\vert \lambda) dt,\;\; \lambda\geq 1,
\end{align*}
where, for $\lambda \geq 1$ and fixed $t$ in support of r.v. $D$,
\begin{align}g_{1,\lambda}(\xi(t),\psi(t),t)&= P_{\boldsymbol{\theta}}[W(\xi(t)Z_1) <W(\psi(t)Z_1)\vert D=t]  \nonumber
	\\&\qquad +\frac{1}{2}P_{\boldsymbol{\theta}}[W(\xi(t)Z_1)\! =\!W(\psi(t)Z_1)\vert D=t].\end{align}
For any fixed $\lambda\geq 1$ and $t$, let $m_{\lambda}^{(1)}(t)$ denote the median of the conditional distributional of $Z_1$ given $D=t$. For any fixed $t$, the conditional p.d.f. of $Z_1$ given $D=t$ is $f_{\lambda}(s\vert t)=\frac{\frac{s}{\lambda}f(s,\frac{st}{\lambda})}{f_D(t\vert \lambda)}$ and $f_D(t\vert \lambda)=\int_{0}^{\infty}\frac{y}{\lambda}f(y,\frac{yt}{\lambda})dy$, $\lambda\geq 1$. Then $\int_{0}^{m_{\lambda}^{(1)}(t)} \, sf(s,\frac{st}{\lambda})ds=\frac{1}{2} \int_{0}^{\infty}sf(s,\frac{st}{\lambda})ds$. It follows from Lemma 4.1 that, for any fixed $t$ and $\lambda\geq 1$, $g_{1,\lambda}(\xi(t),\psi(t),t)> \frac{1}{2}$, provided $m_{\lambda}^{(1)}(t)\leq \frac{1}{\xi(t)}<\frac{1}{\psi(t)}$ or $ \frac{1}{\psi(t)}<\frac{1}{\xi(t)}\leq m_{\lambda}^{(1)}(t)$. Also, for any fixed $t$, $g_{1,\lambda}(\psi(t),\psi(t),t)= \frac{1}{2}, \; \forall \; \lambda\geq 1$.\vspace*{2mm}

On similar lines as in Theorem 3.1.1, under certain conditions, the following theorem provides shrinkage type improvements over an arbitrary scale equivariant estimator under the GPN criterion with a general loss function.\vspace*{2mm}

\noindent
\textbf{Theorem 4.1.1.} Let $\delta_{\psi}(\bold{X})=\psi(D)X_1$ be a scale equivariant estimator of $\theta_1$. Let $l^{(1)}(t)$ and $u^{(1)}(t)$ be functions such that $0<l^{(1)}(t)\leq m_{\lambda}^{(1)}(t)\leq u^{(1)}(t), \;\forall\;\lambda\geq 1$ and any $t$. For any fixed $t$, define $\psi^{*}(t)\!=\!\max\{\frac{1}{u^{(1)}(t)},\min\{\psi(t),\frac{1}{l^{(1)}(t)}\}\}$. Then, under the GPN criterion,
the estimator $\delta_{\psi^{*}}(\bold{X})\!=\!\psi^{*}(D)X_1$ is nearer to $\theta_1$ than the estimator $\delta_{\psi}(\bold{X})=\psi(D)X_1$, provided $P_{\boldsymbol{\theta}}\left[\frac{1}{u^{(1)}(D)}\leq \psi(D)\leq \frac{1}{l^{(1)}(D)}\right]<1,\;\forall \; \boldsymbol{\theta}\in\Theta_0$.
\begin{proof}\!\!\!.
	The GPN of the estimator $\delta_{\psi^{*}}(\bold{X})=\psi^{*}(D)X_1$ relative to $\delta_{\psi}(\bold{X})=\psi(D)X_1$ can be written as 
	\begin{align*}
		GPN(\delta_{\psi^{*}},\delta_{\psi};\boldsymbol{\theta})=\int_{0}^{\infty} g_{1,\lambda}(\psi^{*}(t),\psi(t),t) f_D(t\vert \lambda) dt,\;\;\lambda\geq 1,
	\end{align*}
	where, for $\lambda \geq 1$ and $t$ in support of r.v. $D$, $g_{1,\lambda}(\cdot,\cdot,\cdot)$ is defined by (4.2).
	
	\noindent
	Let $A=\{t:\psi(t)<\frac{1}{u^{(1)}(t)}\}$, $B=\{t:\frac{1}{u^{(1)}(t)}\leq \psi(t)\leq \frac{1}{l^{(1)}(t)}\}$ and $C=\{t:\psi(t)>\frac{1}{l^{(1)}(t)}\}$. Then 
	
	$$\psi^{*}(t)=\begin{cases} \frac{1}{u^{(1)}(t)}, & t\in A\\  \psi(t), & t\in B\\ \frac{1}{l^{(1)}(t)}, & t\in C\end{cases}.$$
	
	Since $l^{(1)}(t)\leq m_{\lambda}^{(1)}(t)\leq u^{(1)}(t)\; (\text{or } \frac{1}{u^{(1)}(t)}\leq \frac{1}{m_{\lambda}^{(1)}(t)}\leq \frac{1}{l^{(1)}(t)}),\; \forall\; \lambda\geq 1$ and $t$, using Lemma 4.1, we have $g_{1,\lambda}(\psi^{*}(t),\psi(t),t)> \frac{1}{2}, \; \forall \; \lambda\geq 1$, whenever $t\in A\cup C$. Also, for $t\in B$, $g_{1,\lambda}(\psi^{*}(t),\psi(t),t)= \frac{1}{2}, \; \forall \; \lambda\geq 1$. Since $P_{\boldsymbol{\theta}}(A\cup C)>0,\; \forall \; \boldsymbol{\theta}\in \Theta_0$, we conclude that
	\\~\\  $GPN(\delta_{\psi^{*}},\delta_{\psi};\boldsymbol{\theta})$
	\begin{align*}
		&\!=\!\int_{A}\! g_{1,\lambda}(\psi^{*}(t),\psi(t),t) f_D(t\vert \lambda) dt \!+ \!\int_{B}\! g_{1,\lambda}(\psi^{*}(t),\psi(t),t) f_D(t\vert \lambda) dt \!\\&\qquad\qquad\qquad\qquad\qquad\qquad\qquad\qquad\qquad\qquad\qquad+\! \int_{C}\! g_{1,\lambda}(\psi^{*}(t),\psi(t),t) f_D(t\vert \lambda) dt\\
		&> \frac{1}{2}, \;\; \boldsymbol{\theta}\in\Theta_0.
	\end{align*}
\end{proof}

\noindent
The proof of the following corollary is contained in the proof of Theorem 4.1.1, and hence skipped.
\\~\\	\textbf{Corollary 4.1.1.} Let $\delta_{\psi}(\bold{X})=\psi(D)X_1$ be a scale equivariant estimator of $\theta_1$. Let $\psi_{1,0}:\Re_{++}\rightarrow \Re_{++}$ be such that $\psi(t)< \psi_{1,0}(t) \leq \frac{1}{u^{(1)}(t)}$, whenever $\psi(t)< \frac{1}{u^{(1)}(t)}$, and $\frac{1}{l^{(1)}(t)}\leq \psi_{1,0}(t)< \psi(t)$, whenever $\frac{1}{l^{(1)}(t)}< \psi(t)$, where $l^{(1)}(\cdot)$ and $u^{(1)}(\cdot)$ are as defined in Theorem 4.1.1. Also let $\psi_{1,0}(t)=\psi(t)$, whenever $\frac{1}{u^{(1)}(t)}\leq \psi(t) \leq \frac{1}{l^{(1)}(t)}$. Then, 
the estimator $\delta_{\psi_{1,0}}(\bold{X})=\psi_{1,0}(D) X_1$ is nearer to $\theta_1$ than $\delta_{\psi}(\bold{X})=\psi(D)X_1$, provided $P_{\boldsymbol{\theta}}\left[\frac{1}{u^{(1)}(D)}\leq \psi(D)\leq \frac{1}{l^{(1)}(D)}\right]<1,\;\forall \; \boldsymbol{\theta}\in\Theta_0$.\vspace*{2mm}

The following corollary provides improvements over the PNSEE $\delta_{1,PNSEE}(\bold{X})=\frac{X_1}{m_{0,1}}$, under the restricted parameter space $\Theta_0$.
\\~\\\textbf{Corollary 4.1.2.} Let $\xi^{*}(t)\!=\!\max\{\frac{1}{u^{(1)}(t)},\min\{\frac{1}{m_{0,1}},\frac{1}{l^{(1)}(t)}\}\}$, where $l^{(1)}(\cdot)$ and $u^{(1)}(\cdot)$ are as defined in Theorem 4.1.1. Then, under the GPN criterion,
the estimator $\delta_{\psi^{*}}(\bold{X})\!=\!\xi^{*}(D)X_1$ is nearer to $\theta_1$ than the PNSEE $\delta_{1,PNSEE}(\bold{X})=\frac{X_1}{m_{0,1}}$, provided $P_{\boldsymbol{\theta}}\left[l^{(1)}(D)\leq m_{0,1}\leq u^{(1)}(D)\right]<1,\;\forall \; \boldsymbol{\theta}\in\Theta_0$.\vspace*{2mm}

In various applications of Theorem 4.1.1 and Corollaries 4.1.1-4.1.2, a common choice for $(l^{(1)}(t),u^{(1)}(t))$ is given by $l^{(1)}(t)=\inf_{\lambda\geq 1} m_{\lambda}^{(1)}(t)$ and $u^{(1)}(t)=\sup_{\lambda\geq 1} m_{\lambda}^{(1)}(t)$.\vspace*{2mm}

In order to identify the behaviour of function $m_{\lambda}^{(1)}(t)$, for any fixed $t$, the following lemma be useful in many situations. Since the proof of the lemma is on the same lines as in Lemma 3.1.1, it is skipped.
\\~\\\textbf{Lemma 4.1.1.} If, for every fixed $\lambda\geq 1$ and $t$, $f(s,\frac{st}{\lambda})/f(s,st)$ is increasing (decreasing) in $s$ (wherever the ratio is not of the form $0/0$), then, for every fixed $t$, $m_{\lambda}^{(1)}(t)$ is an increasing (decreasing) function of $\lambda\in [1,\infty)$.\vspace*{2mm}

\noindent
Under the assumptions of Lemma 4.1.1, one may take, for any fixed $t$, \begin{align}\label{eq:4.2}
	l^{(1)}(t)& =\inf_{\lambda\geq 1} m_{\lambda}^{(1)}(t) =m_1^{(1)}(t)\; (=\lim_{\lambda\to \infty} m_{\lambda}^{(1)}(t))\\
	\;\text{and}\;\;  \label{eq:4.3}
	u^{(1)}(t)&=\sup_{\lambda\geq 1} m_{\lambda}^{(1)}(t)=\lim_{\lambda\to \infty} m_{\lambda}^{(1)}(t)\;(=m_1^{(1)}(t)),
\end{align}
while applying Theorem 4.1.1 and Corollary 4.1.1. \vspace*{2mm}

\vspace*{2mm}

Now we will consider some applications of Theorem 4.1.1 and Corollaries 4.1.1-4.1.2 to specific probability models.

\vspace*{2.5mm}

\noindent
\textbf{Example 4.1.1.} Let $X_1$ and $X_2$ be independent gamma random variables with joint p.d.f. (\ref{eq:3.1}), where, $f(z_1,z_2)=\frac{z_1^{\alpha_1-1}z_2^{\alpha_2-1}e^{-z_1}e^{-z_2}}{\Gamma(\alpha_1)\Gamma(\alpha_2)},\; (z_1,z_2)\in\Re_{++}^2$, for known positive constants $\alpha_1$ and $\alpha_2.$\vspace*{1.5mm}

Consider estimation of the smaller scale parameter $\theta_1$, under the GPN criterion with a general loss function, when it is known apriori that $\boldsymbol{\theta}\in \Theta_0$. Here the restricted MLE of $\theta_1$ is $\delta_{1,RMLE}(\bold{X})\!=\!\min\big\{\!\frac{X_1}{\alpha_1},\frac{X_1+X_2}{\alpha_1+\alpha_2}\!\big\}=X_1\psi_{1,RMLE}(D)$, where $\psi_{1,RMLE}(D)=\min\big\{\frac{1}{\alpha_1},\frac{1+D}{\alpha_1+\alpha_2}\}$ and the unrestricted PNSEE of $\theta_1 \text{ is } \delta_{1,PNSEE}(\bold{X})=\frac{X_1}{m_{0,1}}$, where $m_{0,1}$ is such that $\frac{1}{\Gamma(\alpha_1)} \int_{0}^{m_{0,1}} t^{\alpha_1-1} e^{-t} \,dt=\frac{1}{2}$.
\\~\\ For $t\in\Re_{++}$ and $\lambda\geq 1$, the conditional p.d.f. of $Z_1$ given $D=t$ is
\begin{align*}
	f_{\lambda}(s\vert t)=\begin{cases}
		\frac{ (1+\frac{t}{\lambda})^{\alpha_1+\alpha_2}\, s^{\alpha_1+\alpha_2-1}\,  e^{-(1+\frac{t}{\lambda})s}}{\Gamma{(\alpha_1+\alpha_2)}}, &\text{ if} \;\; 0<s<\infty \\
		0, & \text{ otherwise}  \end{cases}.
\end{align*}
For $t\in (0,\infty)$ and $\lambda\geq 1$, let $m_{\lambda}^{(1)}(t)$ be the median of the p.d.f. $f_{\lambda}(s\vert t),\; t\in (0,\infty)$ and $\lambda\geq 1.$ 
For $\alpha>0$, let $\nu(\alpha)$ denote the median of Gamma($\alpha$,1) distribution, i.e. $\frac{1}{\Gamma(\alpha)} \int_{0}^{\nu(\alpha)} t^{\alpha-1} e^{-t} \,dt=\frac{1}{2}$. Then, $m_{0,1}=\nu(\alpha_1)$ and, for any $t>0$ and $\lambda\geq 1$, $\left(\frac{\lambda+t}{\lambda}\right) m_{\lambda}^{(1)}(t)=\nu(\alpha_1+\alpha_2)$, $m_1^{(1)}(t)=\frac{\nu(\alpha_1+\alpha_2)}{1+t}$ and $\lim_{\lambda\to \infty} m_{\lambda}(t)=\nu(\alpha_1+\alpha_2)$. From Chen and Rubin (\citeyear{MR858317}), we have $\alpha_1+\alpha_2-\frac{1}{3}<\nu(\alpha_1+\alpha_2)<\alpha_1+\alpha_2$.

Thus, as in \eqref{eq:4.2} and \eqref{eq:4.3}, we may take
$l^{(1)}(t)=	m_{1}^{(1)}(t)=\frac{\nu(\alpha_1+\alpha_2)}{1+t}
\text{ and }
u^{(1)}(t)\! 
=\lim\limits_{\lambda\to \infty}	m_{\lambda}^{(1)}(t)=\nu(\alpha_1+\alpha_2).$
\vspace*{1.5mm}

%	Let $\delta_{\psi}(\bold{X})=X_1\psi(D)$ be any scale equivariant estimator of $\theta_1$. Define 
%	$\xi^{*}(t)\!=\!\max\{\frac{1}{u^{(1)}(t)},\min\{\psi(t),\frac{1}{l^{(1)}(t)}\}\},\; t\in \Re_{++}$.
%	By an application of Theorem 4.1.1, the estimator $\delta_{\xi^{*}}(\bold{X})=\xi^{*}(D)X_1$ is better than the equivariant estimator $\delta_{\psi}(\bold{X})=\psi(D) X_1$ under the GPN criterion. \vspace*{2mm}
The following conclusions are immediate from Theorem 4.1.1 and Corollaries 4.1.1-4.1.2:	
\\~\\ \textbf{(i)}  The estimator $\delta_{1,RMLE}^{*}(\bold{X})=\max\big\{\frac{X_1}{\nu(\alpha_1+\alpha_2)},\min\big\{\frac{X_1}{\alpha_1}, \frac{X_1+X_2}{\alpha_1+\alpha_2}\big\}\big\}$ is nearer to $\theta_1$ than the restricted MLE $\delta_{1,RMLE}(\bold{X})=\min\{\frac{X_1}{\alpha_1}, \frac{X_1+X_2}{\alpha_1+\alpha_2}\}$. Clearly, for $\nu(\alpha_1+\alpha_2)\geq \alpha_1$,
$$\delta_{1,RMLE}^{*}(\bold{X})=\begin{cases} \frac{X_1}{\nu(\alpha_1+\alpha_2)}, &\text{ if  } 0<\frac{X_2}{X_1}\leq \frac{\alpha_1+\alpha_2}{\nu(\alpha_1+\alpha_2)}-1\\ 
	\frac{X_1+X_2}{\alpha_1+\alpha_2}, &\text{ if  } \frac{\alpha_1+\alpha_2}{\nu(\alpha_1+\alpha_2)}-1<\frac{X_2}{X_1}\leq \frac{\alpha_2}{\alpha_1}\\
	\frac{X_1}{\alpha_1}, &\text{ if  } \frac{X_2}{X_1}> \frac{\alpha_2}{\alpha_1} \end{cases} $$
and, for $\nu(\alpha_1+\alpha_1)<\alpha_1$, $\delta_{1,RMLE}^{*}(\bold{X})= \frac{X_1}{\nu(\alpha_1+\alpha_2)}$. Since $\nu(\alpha_1+\alpha_2)>\alpha_1+\alpha_2-\frac{1}{3}$, we have $\nu(\alpha_1+\alpha_2)> \alpha_1$, if $\alpha_2\geq \frac{1}{3}$.
\\~\\ \textbf{(ii)}   The estimator $\delta_{1,PNSEE}^{*}(\bold{X})= \max\big\{\frac{X_1}{\nu(\alpha_1+\alpha_2)},\min\big\{\frac{X_1}{\nu(\alpha_1)}, \frac{X_1+X_2}{\nu(\alpha_1+\alpha_2)}\big\}\big\}$ is nearer to $\theta_1$ than $\delta_{1,PNSEE}(\bold{X})=\frac{X_1}{\nu(\alpha_1)}$. Clearly
$$\delta_{1,PNSEE}^{*}(\bold{X})=\begin{cases} 
	\frac{X_1+X_2}{\nu(\alpha_1+\alpha_2)}, &\text{ if  } 0<\frac{X_2}{X_1}\leq \frac{\nu(\alpha_1+\alpha_2)}{\nu(\alpha_1)}-1\\
	\frac{X_1}{\nu(\alpha_1)}, &\text{ if  } \frac{X_2}{X_1}>\frac{\nu(\alpha_1+\alpha_2)}{\nu(\alpha_1)}-1 \end{cases}. $$
\textbf{(iii)}   The estimator $\delta_{1,UE}^{*}(\bold{X})= \max\big\{\frac{X_1}{\nu(\alpha_1+\alpha_2)},\min\big\{\frac{X_1}{\alpha_1}, \frac{X_1+X_2}{\nu(\alpha_1+\alpha_2)}\big\}\big\}$ is nearer to $\theta_1$ than the unbiased estimator $\delta_{1,UE}(\bold{X})=\frac{X_1}{\alpha_1}$. Clearly, for $\nu(\alpha_1+\alpha_2)\geq \alpha_1$,
$$\delta_{1,UE}^{*}(\bold{X})=\begin{cases} \frac{X_1+X_2}{\nu(\alpha_1+\alpha_2)}, &\text{ if  } 0\leq \frac{X_2}{X_1}< \frac{\nu(\alpha_1+\alpha_2)}{\alpha_1}-1\\ 
	\frac{X_1}{\alpha_1}, &\text{ if  } \frac{X_2}{X_1}\geq \frac{\nu(\alpha_1+\alpha_2)}{\alpha_1}-1 \end{cases} $$
and, for $\nu(\alpha_1+\alpha_1)<\alpha_1$, $\delta_{1,UE}^{*}(\bold{X})= \frac{X_1}{\nu(\alpha_1+\alpha_2)}$. \vspace*{2mm}

\noindent
\textbf{(iv)}  For $\nu(\alpha_1+\alpha_2)\leq \alpha_1$, the estimator $\delta_{1,UE}(\bold{X})= \frac{X_1}{\alpha_1}$ is nearer to $\theta_1$ than the restricted MLE $\delta_{1,RMLE}(\bold{X})=\min\{\frac{X_1}{\alpha_1}, \frac{X_1+X_2}{\alpha_1+\alpha_2}\}$. Ma and Liu (2014) proved a similar result under the GPN criterion with a specific loss function $L_1(\boldsymbol{\theta},a)=\vert \frac{a}{\theta_1}-1\vert,\; \boldsymbol{\theta}\in \Theta_0,\; a>0$. In fact several results reported in Ma and Liu (2014) can be obtained as particular cases of Corollary 4.1.1.

\vspace*{2mm}

\noindent
\textbf{Example 4.1.2.} Let $X_1$ and $X_2$ be independent random variables with joint p.d.f. (\ref{eq:3.1}),
where $\boldsymbol{\theta}=(\theta_1,\theta_2)\in \Theta_0$ and $f\left(z_1,z_2\right)= \alpha_1 \alpha_2 z_1^{\alpha_1-1}z_2^{\alpha_2-1},\text{ if } 0<z_1 <1,\,0<z_2 <1; =0, \text{ otherwise}$, for known positive constants $\alpha_1$ and $\alpha_2$.
\vspace*{1.5mm}

Consider estimation of $\theta_1$ under the GPN criterion with general loss function. Here the unrestricted PNSEE of $\theta_1 \text{ is } \delta_{1,PNSEE}(\bold{X})=2^{\frac{1}{\alpha_1}}X_1$. Also, for $t\in\Re_{++}$, the conditional p.d.f. of $Z_1$ given $D=t$ is
%	\newpage
\begin{align*}
	f_{\lambda}(s\vert t)=\frac{(\alpha_1+\alpha_2)s^{\alpha_1+\alpha_2-1}}{\left(\min\{1,\frac{\lambda}{t}\}\right)^{\alpha_1+\alpha_2}}, \text{ if} \; 0<s<\min\bigg\{1,\frac{\lambda}{t}\bigg\},\; \lambda \geq 1.
\end{align*}

The median of the above density is  $m_{\lambda}^{(1)}(t)=2^{\frac{-1}{\alpha_1+\alpha_2}}\min\big\{1,\frac{\lambda}{t}\big\},\;\forall \; t\in \Re_{++}$ and $\lambda\geq 1.$ Clearly, for $t\in\Re_{++}$, $m_{\lambda}^{(1)}(t)$ is increasing in $\lambda\in[1,\infty)$.
As in \eqref{eq:4.2} and \eqref{eq:4.3}, we take
$l^{(1)}(t)=	m_{1}^{(1)}(t)= 2^{\frac{-1}{\alpha_1+\alpha_2}}\min\big\{1,\frac{1}{t}\big\}
\text{ and }
u^{(1)}(t)\! 
=\lim\limits_{\lambda\to \infty}	m_{\lambda}^{(1)}(t)=2^{\frac{-1}{\alpha_1+\alpha_2}}.$
\vspace*{1.5mm}

%Let $\delta_{\psi}(\bold{X})=X_1\psi(D)$ be any scale equivariant estimator of $\theta_1$. Then, using Theorem 4.1.1, $GPN(\delta_{\xi^{*}},\delta_{\psi};\boldsymbol{\theta})\geq \frac{1}{2},\; \forall \; \boldsymbol{\theta}\in \Theta_0$, where $\delta_{\xi^{*}}(\bold{X})=\xi^{*}(D)X_1$ and
%$\xi^{*}(t)=\max\!\big\{\!2^{\frac{1}{\alpha_1+\alpha_2}},\min\{2^{\frac{1}{\alpha_1+\alpha_2}}\max\{1,t\},\psi(t)\}\big\},\; t>0.$\vspace*{2.5mm}

\noindent The following conclusions immediately follow from Theorem 4.1.1 and Corollary 4.1.1:
\\~\\\textbf{(i)}  Define $\delta_{1,PNSEE}^{*}(\bold{X})=\max\!\big\{\!2^{\frac{1}{\alpha_1+\alpha_2}},\min\{2^{\frac{1}{\alpha_1+\alpha_2}}\max\{1,t\},2^{\frac{1}{\alpha_1}}\}\big\} X_1=$ $\min\{2^{\frac{1}{\alpha_1+\alpha_2}}\max\{1,t\},2^{\frac{1}{\alpha_1}}\} X_1$. Then the estimator $\delta_{1,PNSEE}^{*}(\bold{X})$ is nearer to $\theta_1$ than the PNSEE $\delta_{1,PNSEE}(\bold{X})$. It is easy to verify that
$$\delta_{1,PNSEE}^{*}(\bold{X})=\begin{cases} 2^{\frac{1}{\alpha_1+\alpha_2}} X_1, &\text{if  } 0<\frac{X_2}{X_1}<1\\2^{\frac{1}{\alpha_1+\alpha_2}} X_2, &\text{if  } 1\leq \frac{X_2}{X_1}<2^{\frac{\alpha_2}{\alpha_1(\alpha_1+\alpha_2)}}\\ 2^{\frac{1}{\alpha_1}} X_1, &\text{if  } \frac{X_2}{X_1}\geq 2^{\frac{\alpha_2}{\alpha_1(\alpha_1+\alpha_2)}} \end{cases}.  $$
\\~\\ \textbf{(ii)}  Let $\psi_{1,0}(t)$ be such that $2^{\frac{1}{\alpha_1+\alpha_2}}\max\{1,t\} \leq \psi_{1,0}(t)< 2^{\frac{1}{\alpha_1}},\; \forall\; t\leq 2^{\frac{\alpha_2}{\alpha_1(\alpha_1+\alpha_2)}}$, and $\psi_{1,0}(t)=2^{\frac{1}{\alpha_1}},\; \forall\; t>2^{\frac{\alpha_2}{\alpha_1(\alpha_1+\alpha_2)}}$. Then the estimator $\delta_{\psi_{1,0}}(\bold{X})=\psi_{1,0}(D)X_1$ is nearer to $\theta_1$ than the PNSEE $\delta_{1,PNSEE}(\bold{X})$.

%Under the scaled mean squared error $\left(\text{i.e., } E[(\frac{\delta}{\theta_1}-1)^2]\right)$ criterion, the unrestricted BSEE is $\delta_{1,c_{0,1}}(\bold{X})=\frac{\alpha_1+2}{\alpha_1+1} X_1$ (i.e., $c_{0,1}=\frac{\alpha_1+2}{\alpha_1+1}$). Using Corollary 4.1.3, we get $GPN(\delta_{1,c_{0,1}}^{*},\delta_{1,c_{0,1}};\boldsymbol{\theta})\geq \frac{1}{2},\; \forall \; \boldsymbol{\theta}\in \Theta_0$, where $\delta_{1,c_{0,1}}^{*}(\bold{X})=\max\!\big\{\!2^{\frac{1}{\alpha_1+\alpha_2}},\min\{2^{\frac{1}{\alpha_1+\alpha_2}}\max\{1,t\},c_{0,1}\}\big\}X_1$.	

\subsection{\textbf{Estimation of The Larger Scale Parameter $\theta_2$}} \label{4.2}

In this section, we consider estimation of the larger scale parameter $\theta_2$ under the GPN criterion with a general loss function $L_2(\boldsymbol{\theta},a)=W(\frac{a}{\theta_2})$, when it is known that $0<\theta_1\leq \theta_2 <\infty$ (i.e., $\boldsymbol{\theta}\in \Theta_0$). Here $W:\Re_{++}\rightarrow \Re_{++}$ is such that $W(1)=0$, $W(t)$ is strictly decreasing in $(0,1)$ and strictly increasing in $(1,\infty)$.
Here, any scale equivariant estimator of $\theta_2$ is of the form
$\delta_{\psi}(\bold{X})=\psi(D)X_2,$ for some function $\psi:\,\Re_{++}\rightarrow \Re_{++}$, where $D=\frac{X_2}{X_1}$.
Define $Z_2=\frac{X_2}{\theta_2}$, $\lambda=\frac{\theta_2}{\theta_1}$ and $f_D(t\vert \lambda)$ the p.d.f. of r.v. $D=\frac{X_2}{X_1}$. Let $\delta_{\xi}(\bold{X})=\xi(D)X_2$ and $\delta_{\psi}(\bold{X})=\psi(D)X_2$ be two scale equivariant estimators of $\theta_2$. Then, the GPN of $\delta_{\xi}(\bold{X})=\xi(D)X_2$ relative to $\delta_{\psi}(\bold{X})=\psi(D)X_2$ can be written as 
\begin{align*}
	GPN(\delta_{\xi},\delta_{\psi};\boldsymbol{\theta})&=\int_{0}^{\infty} g_{2,\lambda}(\xi(t),\psi(t),t) f_D(t\vert \lambda) dt,\;\; \boldsymbol{\theta}\in \Theta_0,
\end{align*}
where, for $\lambda \geq 1$, $g_{2,\lambda}(\xi(t),\psi(t),t)= P_{\boldsymbol{\theta}}[W(\xi(t)Z_2) <W(\psi(t)Z_2)\vert D=t]
+\frac{1}{2}P_{\boldsymbol{\theta}}[W(\xi(t)Z_2)$  $ = W(\psi(t)Z_2)\vert D=t].$
For any fixed $\lambda\geq 1$ and $t$, let $m_{\lambda}^{(2)}(t)$ denote the median of the conditional distributional of $Z_2$ given $D=t$. For any fixed $t$, the conditional p.d.f. of $Z_2$ given $D=t$ is $f_{\lambda}(s\vert t)=\frac{\frac{\lambda s}{t^2}f(\frac{\lambda s}{t},s)}{f_D(t\vert \lambda)}$ and $f_D(t\vert \lambda)=\int_{0}^{\infty}\frac{\lambda y}{t^2}f(\frac{\lambda y}{t},y)dy$, $\lambda\geq 1$. Thus $\int_{0}^{m_{\lambda}^{(2)}(t)} sf(\frac{\lambda s}{t},s)ds=\frac{1}{2} \int_{0}^{\infty}sf(\frac{\lambda s}{t},s)ds$. Using Lemma 4.1, we have, for any fixed $t$ and $\lambda\geq 1$, $g_{2,\lambda}(\xi(t),\psi(t),t)> \frac{1}{2}$, provided $m_{\lambda}^{(2)}(t)\leq \frac{1}{\xi(t)}<\frac{1}{\psi(t)}$ or $ \frac{1}{\psi(t)}<\frac{1}{\xi(t)}\leq m_{\lambda}^{(2)}(t)$. Moreover, for any fixed $t$, $g_{2,\lambda}(\psi(t),\psi(t),t)= \frac{1}{2}, \; \forall \; \lambda\geq 1$. These arguments lead to the following results.

\vspace*{2mm}

\noindent
\textbf{Theorem 4.2.1.} Let $\delta_{\psi}(\bold{X})=\psi(D)X_2$ be a scale equivariant estimator of $\theta_2$. Let $l^{(2)}(t)$ and $u^{(2)}(t)$ be functions such that $0<l^{(2)}(t)\leq m_{\lambda}^{(2)}(t)\leq u^{(2)}(t), \;\forall\;\lambda\geq 1$ and any $t$. For any fixed $t$, define $\psi^{*}(t)\!=\!\max\{\frac{1}{u^{(2)}(t)},\min\{\psi(t),\frac{1}{l^{(2)}(t)}\}\}$. Then, under the GPN criterion,
the estimator $\delta_{\psi^{*}}(\bold{X})\!=\!\psi^{*}(D)X_2$ is nearer to $\theta_2$ than the estimator $\delta_{\psi}(\bold{X})=\psi(D)X_2$, provided $P_{\boldsymbol{\theta}}\left[\frac{1}{u^{(2)}(D)}\leq \psi(D)\leq \frac{1}{l^{(2)}(D)}\right]<1,\;\forall \; \boldsymbol{\theta}\in\Theta_0$.\vspace*{2mm}

\noindent
\textbf{Corollary 4.2.1.} Let $\delta_{\psi}(\bold{X})=\psi(D)X_2$ be a scale equivariant estimator of $\theta_2$. Let $\psi_{2,0}:\Re_{++}\rightarrow \Re_{++}$ be such that $\psi(t)< \psi_{2,0}(t) \leq \frac{1}{u^{(2)}(t)}$, whenever $\psi(t)< \frac{1}{u^{(2)}(t)}$, and $\frac{1}{l^{(2)}(t)}\leq \psi_{2,0}(t)< \psi(t)$, whenever $\frac{1}{l^{(2)}(t)}< \psi(t)$, where $l^{(2)}(\cdot)$ and $u^{(2)}(\cdot)$ are as defined in Theorem 4.2.1. Also let $\psi_{2,0}(t)=\psi(t)$, whenever $\frac{1}{u^{(2)}(t)}\leq \psi(t) \leq \frac{1}{l^{(2)}(t)}$. Then, 
$GPN(\delta_{\psi_{2,0}},\delta_{\psi};\boldsymbol{\theta})> \frac{1}{2},\; \forall \; \boldsymbol{\theta}\in \Theta_0,$ provided $P_{\boldsymbol{\theta}}\left[\frac{1}{u^{(2)}(D)}\leq \psi(D)\leq \frac{1}{l^{(2)}(D)}\right]<1,\;\forall \; \boldsymbol{\theta}\in\Theta_0$,
where $\delta_{\psi_{2,0}}(\bold{X})=\psi_{2,0}(D)X_2$.\vspace*{2mm}

Note that, in the unrestricted case $\Theta\!=\!\Re_{++}^2$, the PNSEE of $\theta_2$ is $\delta_{2,PNSEE}(\bold{X})=\frac{X_2}{m_{0,2}}$, where $m_{0,2}>0$ is the median of the r.v. $Z_2=\frac{X_2}{\theta_2}$.
The following corollary provides improvements over the PNSEE under the restricted parameter space.
\\~\\\textbf{Corollary 4.2.2.} Let $\xi^{*}(t)\!=\!\max\{\frac{1}{u^{(2)}(t)},\min\{\frac{1}{m_{0,2}},\frac{1}{l^{(2)}(t)}\}\}$, where $l^{(2)}(\cdot)$ and $u^{(2)}(\cdot)$ are as defined in Theorem 4.2.1. Then, under the GPN criterion,
the estimator $\delta_{\xi^{*}}(\bold{X})\!=\!\xi^{*}(D)X_2$ is nearer to $\theta_2$ than the PNSEE $\delta_{2,PNSEE}(\bold{X})=\frac{X_2}{m_{0,2}}$, provided $P_{\boldsymbol{\theta}}\left[l^{(2)}(D)\leq m_{0,2}\leq u^{(2)}(D)\right]<1,\;\forall \; \boldsymbol{\theta}\in\Theta_0$.\vspace*{2mm}

In order to identify the behaviour of $m_{\lambda}^{(2)}(t),$ for any fixed $t$, we have the following lemma on the lines of Lemma 4.1.1.
\\~\\\textbf{Lemma 4.2.1.} If, for every fixed $\lambda\geq 1$ and $t$, $f(\frac{\lambda s}{t},s)/f(\frac{ s}{t},s)$ is increasing (decreasing) in $s$ (wherever the ratio is not of the form $0/0$), then, for every fixed $t$, $m_{\lambda}^{(2)}(t)$ is an increasing (decreasing) function of $\lambda\in [1,\infty)$.\vspace*{2mm}

\noindent
Under the assumptions of Lemma 4.2.1, one may take, for any fixed $t$, \begin{align}\label{eq:4.4}
	l^{(2)}(t)& =\inf_{\lambda\geq 1} m_{\lambda}^{(2)}(t) =m_1^{(2)}(t)\; (=\lim_{\lambda\to \infty} m_{\lambda}^{(2)}(t))\\
	\;\text{and}\;\; \label{eq:4.5}
	u^{(2)}(t)&=\sup_{\lambda\geq 1} m_{\lambda}^{(2)}(t)=\lim_{\lambda\to \infty} m_{\lambda}^{(2)}(t)\;(=m_1^{(2)}(t)),
\end{align}
while applying Theorem 4.2.1 and Corollary 4.2.1. 
%\vspace*{2mm}

As in Section 4.1, we will now apply Theorem 4.2.1 and Corollaries 4.2.1-4.2.2 to estimation of the larger scale parameter $\theta_2$ in scale probability models considered in Examples 4.1.1-4.1.2.

\vspace*{2mm}

\noindent \textbf{Example 4.2.1.} Let $X_1$ and $X_2$ be independent gamma random variables as defined in Example 4.1.1. 
Consider estimation of $\theta_2$, under the GPN criterion with a general loss function. 
\\~\\ For $t\in\Re_{++}$ and $\lambda\geq 1$, the conditional p.d.f. of $Z_2$ given $D=t$ is
\begin{align*}
	f_{\lambda}(s\vert t)=\begin{cases}
		\frac{ (1+\frac{\lambda}{t})^{\alpha_1+\alpha_2}\, s^{\alpha_1+\alpha_2-1}\,  e^{-(1+\frac{\lambda}{t})s}}{\Gamma{(\alpha_1+\alpha_2)}}, &\text{ if} \;\; 0<s<\infty \\
		0, & \text{ otherwise}  \end{cases}.
\end{align*}
Let $\nu(\alpha)$ denote the median of Gamma($\alpha$,1) distribution. Then $ m_{\lambda}^{(2)}(t)=\left(\frac{t}{\lambda+t}\right)\nu(\alpha_1+\alpha_2)$, $\lambda\geq 1,\; t>0$, and, as in (4.5) and (4.6), we may take $l^{(2)}(t)=0,\; t>0$ and $u^{(2)}(t)=\frac{t}{1+t}\,\nu(\alpha_1+\alpha_2),\; t>0$. Also $m_{0,2}=\nu(\alpha_2)$ and the PNSEE of $\theta_2$ is $\delta_{2,PNSEE}(\bold{X})=\frac{X_2}{\nu(\alpha_2)}$. The restricted MLE of $\theta_2$ is $\delta_{2,RMLE}(\bold{X})=\max\{\frac{X_2}{\alpha_2}, \frac{X_1+X_2}{\alpha_1+\alpha_2}\}=\psi_{2,R}(D)X_2$, where $\psi_{2,R}(t)=\max\{\frac{1}{\alpha_2}, \frac{1+t}{t(\alpha_1+\alpha_2)}\},\; t>0$. Using Theorem 4.2.2 and Corollary 4.2.1, the following conclusions are evident:	
\\~\\ \textbf{(i)}  The estimator $\delta_{2,RMLE}^{*}(\bold{X})=\max\big\{\frac{X_2}{\alpha_2},\frac{X_1+X_2}{\nu(\alpha_1+\alpha_2)}\big\}$ is nearer to $\theta_2$ than the restricted MLE $\delta_{2,RMLE}(\bold{X})=\max\{\frac{X_2}{\alpha_2}, \frac{X_1+X_2}{\alpha_1+\alpha_2}\}$.
\\~\\ \textbf{(ii)}   The estimator $\delta_{2,PNSEE}^{*}(\bold{X})= \max\big\{\frac{X_2}{\nu(\alpha_2)},\frac{X_1+X_2}{\nu(\alpha_1+\alpha_2)}\big\}$ is nearer to $\theta_2$ than the PNSEE $\delta_{2,PNSEE}(\bold{X})=\frac{X_2}{\nu(\alpha_2)}$. 
\\~\\\textbf{(iii)}   The restricted MLE $\delta_{2,RMLE}(\bold{X})=\max\{\frac{X_2}{\alpha_2}, \frac{X_1+X_2}{\alpha_1+\alpha_2}\}$ is nearer to $\theta_2$ than the unbiased estimator $\delta_{2,UE}(\bold{X})=\frac{X_2}{\alpha_2}$. This result, under a specific loss function $L_2(\boldsymbol{\theta},a)=\vert \frac{a}{\theta_2}-1\vert,\; \boldsymbol{\theta}\in \Theta_0,\; a>0$, is proved in Ma and Liu (\citeyear{ma2014}).

\vspace*{2mm}

\noindent \textbf{Example 4.2.2.} Let $X_1$ and $X_2$ be independent random variables as described in Example 4.1.2. Consider estimation of $\theta_2$ under the GPN criterion with a general error loss function. Here the PNSEE of $\theta_2 \text{ is } \delta_{2,PNSEE}(\bold{X})=2^{\frac{1}{\alpha_2}}X_2$. Also, for $\lambda\geq 1$ and $t\in [0,\infty)$, the conditional p.d.f. of $Z_2$ given $D=t$ is
\begin{align*}
	f_{\lambda}(s\vert t)=\frac{(\alpha_1+\alpha_2)s^{\alpha_1+\alpha_2-1}}{\left(\min\{1,\frac{t}{\lambda}\}\right)^{\alpha_1+\alpha_2}}, \text{ if} \; 0<s<\min\bigg\{1,\frac{t}{\lambda}\bigg\}.
\end{align*}
The median of the above density is  $m_{\lambda}^{(2)}(t)=2^{\frac{-1}{\alpha_1+\alpha_2}}\min\big\{1,\frac{t}{\lambda}\big\},\;\forall \; t\in \Re_{++}$ and $\lambda\geq 1.$ Clearly, for $t\in\Re_{++}$, $m_{\lambda}^{(2)}(t)$ is decreasing in $\lambda\in[1,\infty)$.
As in \eqref{eq:4.4} and \eqref{eq:4.5}, we may take
$l^{(2)}(t)=\lim\limits_{\lambda\to \infty}	m_{\lambda}^{(2)}(t)= 0
\text{ and }
u^{(2)}(t)\! 
=	m_{1}^{(2)}(t)=2^{\frac{-1}{\alpha_1+\alpha_2}}\min\{1,t\}.$\vspace*{1.5mm}

%Let $\delta_{\psi}(\bold{X})=X_2\psi(D)$ be any scale equivariant estimator of $\theta_2$. Define $\delta_{\xi^{*}}(\bold{X})=\xi^{*}(D)X_2$, where
%$\xi^{*}(t)=\max\!\big\{\!2^{\frac{1}{\alpha_1+\alpha_2}}\max\big\{1,\frac{1}{t}\big\},\psi(t)\big\},\; t>0.$ Then, using Theorem 4.3.1, $GPN(\delta_{\xi^{*}},\delta_{\psi};\boldsymbol{\theta})\geq \frac{1}{2},\; \forall \; \boldsymbol{\theta}\in \Theta_0$.\vspace*{2.5mm}

\vspace*{1.5mm}

\noindent The following conclusions immediately follow from Theorem 4.2.1 and Corollary 4.2.1:
\\~\\\textbf{(i)} The estimator $\delta_{2,PNSEE}^{*}(\bold{X})=\max\!\big\{\!2^{\frac{1}{\alpha_1+\alpha_2}}X_1,2^{\frac{1}{\alpha_2}}X_2\big\}$ is nearer to $\theta_1$ than the PNSEE $\delta_{2,PNSEE}(\bold{X})=2^{\frac{1}{\alpha_2}}X_2$.
\\~\\ \textbf{(ii)}  Let $\psi_{2,0}(t)$ be such that $2^{\frac{1}{\alpha_2}} < \psi_{2,0}(t)\leq 2^{\frac{1}{\alpha_1+\alpha_2}}\max\big\{1,\frac{1}{t}\big\},\; \forall\; t\leq 2^{\frac{-\alpha_1}{\alpha_2(\alpha_1+\alpha_2)}}$, and $\psi_{2,0}(t)=2^{\frac{1}{\alpha_2}},\; \forall\; t>2^{\frac{-\alpha_1}{\alpha_2(\alpha_1+\alpha_2)}}$. Then the estimator $\delta_{\psi_{2,0}}(\bold{X})=\psi_{2,0}(D)X_2$ is nearer to $\theta_2$ than the PNSEE $\delta_{2,PNSEE}(\bold{X})=2^{\frac{1}{\alpha_2}}X_2$.

%Under the scaled mean squared error $\left(\text{i.e., } E[(\frac{\delta}{\theta_2}-1)^2]\right)$ criterion, the unrestricted BSEE is $\delta_{2,c_{0,2}}(\bold{X})=\frac{\alpha_2+2}{\alpha_2+1} X_2$ (i.e., $c_{0,2}=\frac{\alpha_2+2}{\alpha_2+1}$). Using Corollary 4.3.3, we get $GPN(\delta_{2,c_{0,2}}^{*},\delta_{2,c_{0,2}};\boldsymbol{\theta})\geq \frac{1}{2},\; \forall \; \boldsymbol{\theta}\in \Theta_0$, where $\delta_{2,c_{0,2}}^{*}(\bold{X})=\max\!\big\{\!2^{\frac{1}{\alpha_1+\alpha_2}}\max\big\{1,\frac{1}{t}\big\},c_{0,2}\big\} X_2$.\vspace*{2.5mm}

\section{\textbf{Simulation Study}}  \label{5}

\subsection{For Smaller Location Parameter $\theta_1$}   \label{5.1}
% Please add the following required packages to your document preamble:
% \usepackage{booktabs}
In Example 3.1.1, we considered a bivariate normal distribution with unknown means $\theta_1$ and $\theta_2$ ($-\infty<\theta_1\leq \theta_2<\infty$), known variances $\sigma_1^2>0$ and $\sigma_2^2>0$, and known correlation coefficient $\rho$ ($-1<\rho<1$). For estimation of $\theta_1$ under GPN criterion, we considered various estimators. To further evaluate the performances of these estimators, in this section, we compare these estimators under the GPN criterion with the absolute error loss (i.e., $W(t)=\vert t\vert,\; t\in \Re$), numerically. For simulations, 10,000 random samples of size 1 were generated from the relevant bivariate normal distribution. For various configurations of $\alpha=\frac{\sigma_2(\sigma_2-\rho\sigma_1)}{\sigma_1^2+\sigma_2^2-2\rho \sigma_1 \sigma_2}$, using the Monte Carlo simulations, we obtained the GPN values of the restricted MLE ($\delta_{1,RMLE}^{*}(\bold{X})=X_1-(1-\alpha) \max\{0,-D\}$) relative to the PNLEE ($\delta_{1,PNLEE}(\bold{X})=X_1$), of the improved Hwang and Peddada (HP) estimator ($\delta_{1,HP}^{*}(\bold{X})=\alpha X_1 +(1-\alpha) X_2$) relative to the Hwang and Peddada (HP) estimator ($\delta_{1,HP}(\bold{X})=X_1-\max\{0,(\alpha-1)D\}$) and of the restricted MLE ($\delta_{1,RMLE}(\bold{X})=X_1-(1-\alpha) \max\{0,-D\}$) relative to the Tan and Peddada (PDT) estimator ($\delta_{1,PDT}(\bold{X})=X_1-\max\{0,-D\}$). These values are tabulated in Tables 1-3.
% After that, we also compared the risk performances, under mean squared error, of the improved restricted MLE (RMLE) ($\delta_{1,R}^{*}$), the improved PNEE ($\delta_{1,m_{0,1}}^{*}$) and improved BLEE ($\delta_{1,c_{0,1}}^{*}$).\vspace*{1.5mm}
%	The simulated values of the risks of improved RMLE ($\delta_{1,R}^{*}$), the improved PNEE ($\delta_{1,m_{0,1}}^{*}$) and improved BLEE ($\delta_{1,c_{0,1}}^{*}$) are plotted in Figure \ref{fig1}. 
%The following observations are evident from Table \ref{table:1}, Table \ref{table:2} and Table \ref{table:3}:
The following observations are evident from Tables 1-3:
\vspace*{2mm}

\noindent (i) All the GPN values are greater than $0.5$, which is in conformity with the theoretical findings of Example 3.1.1.
\\ (ii)  From Table \ref{table:1}, we can observe that, when $\sigma_1$ is relatively larger than $\sigma_2$, we get higher GPN values. Also, for negative $\rho$, the GPN values are higher. 
\\ (iii)  From Table \ref{table:2}, we can see that, as the value of $\rho \sigma_2-\sigma_1$ ($>0$) increases, the GPN value also increases. Also, from Table \ref{table:3}, we observe that as the value of $\rho \sigma_1-\sigma_2$ ($>0$) increases, the GPN value also increases.
%	\newpage
\begin{table}[h!]
	\centering
	\caption{The GPN values of the restricted MLE ($\delta_{1,RMLE}$) relative to the PNLEE ($\delta_{1,PNLEE}$)}
	\label{table:1}
	\centering
	\begin{tabular}{@{}ccccccc@{}}
		\toprule
		& \multicolumn{6}{c}{($\sigma_1,\sigma_2,\rho$)}            \\
		$\theta_2-\theta_1$ & (3,0.5,-0.9)    & (0.5,5,-0.5)    & (1,1,0)    & (15,2,0.2)   & (1,30,0.5)   & (30,1,0.9)   \\          \midrule
		0.0 & 0.743 & 0.557 & 0.560 & 0.708 & 0.549 & 0.740 \\
		0.5 & 0.713 & 0.577 & 0.609 & 0.718 & 0.537 & 0.749 \\
		1.0 & 0.693 & 0.548 & 0.580 & 0.722 & 0.545 & 0.740 \\
		1.5 & 0.662 & 0.558 & 0.547 & 0.719 & 0.535 & 0.740 \\
		2.0 & 0.626 & 0.566 & 0.540 & 0.717 & 0.546 & 0.753 \\
		2.5 & 0.611 & 0.570 & 0.516 & 0.721 & 0.557 & 0.741 \\
		3.0 & 0.584 & 0.575 & 0.507 & 0.709 & 0.556 & 0.732 \\ \bottomrule
	\end{tabular}
	
\end{table}

% Please add the following required packages to your document preamble:
% \usepackage{booktabs}
\begin{table}[h!]
	\caption{The GPN values of the improved HP ($\delta_{1,HP}^{*}$) relative to the HP ($\delta_{1,HP}$) when $\alpha>1$ (i.e. $\sigma_1<\rho \sigma_2$)} 
	\label{table:2}
	\centering
	\begin{tabular}{@{}ccccccc@{}}
		\toprule
		& \multicolumn{6}{c}{($\sigma_1,\sigma_2,\rho$)}              \\
		$\theta_2-\theta_1$ & (0.1,5,0.2)    & (1,25,0.2)   & (0.5,2,0.5)   & (5,15,0.5)   & (0.5,5,0.9)   & (2,15,0.9)   \\  \midrule
		0   & 0.520 & 0.502 & 0.524 & 0.512 & 0.620 & 0.619 \\
		0.5 & 0.515 & 0.514 & 0.523 & 0.521 & 0.631 & 0.620 \\
		1   & 0.516 & 0.508 & 0.532 & 0.518 & 0.633 & 0.621 \\
		1.5 & 0.520 & 0.502 & 0.528 & 0.520 & 0.639 & 0.625 \\
		2   & 0.524 & 0.511 & 0.521 & 0.514 & 0.626 & 0.634 \\
		2.5 & 0.520 & 0.520 & 0.518 & 0.513 & 0.621 & 0.630 \\
		3   & 0.513 & 0.515 & 0.510 & 0.516 & 0.609 & 0.629  \\ \bottomrule
	\end{tabular}
\end{table}

\begin{table}[h!]
	\caption{The GPN values of the restricted MLE ($\delta_{1,RMLE}$) relative to the PDT ($\delta_{1,PDT}$) when $\alpha <0$ (i.e. $\sigma_2<\rho \sigma_1$)}
	\label{table:3}
	\centering
	\begin{tabular}{@{}ccccccc@{}}
		\toprule
		& \multicolumn{6}{c}{($\sigma_1,\sigma_2,\rho$)}              \\
		$\theta_2-\theta_1$ & (5,0.1,0.2)    & (25,1,0.2)   & (2,0.5,0.5)   & (15,5,0.5)   & (5,0.5,0.9)   & (15,2,0.9)   \\  \midrule
		0   & 0.512 & 0.516 & 0.522 & 0.516 & 0.617 & 0.619 \\
		0.5 & 0.733 & 0.613 & 0.650 & 0.526 & 0.723 & 0.679 \\
		1   & 0.706 & 0.677 & 0.640 & 0.550 & 0.708 & 0.712 \\
		1.5 & 0.692 & 0.713 & 0.598 & 0.577 & 0.683 & 0.722 \\
		2   & 0.672 & 0.729 & 0.569 & 0.588 & 0.667 & 0.720 \\
		2.5 & 0.653 & 0.730 & 0.539 & 0.597 & 0.644 & 0.710 \\
		3   & 0.633 & 0.725 & 0.522 & 0.611 & 0.630 & 0.706 \\   \bottomrule
	\end{tabular}
\end{table}

\subsection{For Larger Scale Parameter $\theta_2$}   \label{5.2}
% Please add the following required packages to your document preamble:
% \usepackage{booktabs}
In this section, for estimation of the larger scale parameter $\theta_2$, under the GPN criterion with loss function $L_2(\underline{\theta},a)=\vert\frac{a}{\theta_2}-1\vert,\; a\in \mathcal{A}=(0,\infty),\; \underline{\theta}\in \Theta_0$, we numerically compare various estimators considered in Example 4.2.1.
For simulations, 10,000 random samples of size 1 were generated from relevant gamma distributions. Using the Monte Carlo simulations, we obtained the GPN values of the improved restricted MLE ($\delta_{2,RMLE}^{*}(\bold{X})=\max\big\{\frac{X_2}{\alpha_2},\frac{X_1+X_2}{\nu(\alpha_1+\alpha_2)}\big\}$) relative to the restricted MLE ($\delta_{2,RMLE}(\bold{X})=\max\big\{\frac{X_2}{\alpha_2},\frac{X_1+X_2}{\alpha_1+\alpha_2}\big\}$), of the improved PNSEE ($\delta_{2,PNSEE}^{*}(\bold{X})= \max\big\{\frac{X_2}{\nu(\alpha_2)},\frac{X_1+X_2}{\nu(\alpha_1+\alpha_2)}\big\}$) relative to the PNSEE ($\delta_{2,PNSEE}(\bold{X})=\frac{X_2}{\nu(\alpha_2)}$) and of the restricted MLE ($\delta_{2,RMLE}(\bold{X})=\max\big\{\frac{X_2}{\alpha_2},\frac{X_1+X_2}{\alpha_1+\alpha_2}\big\}$) relative to the unbiased estimator ($\delta_{2,UE}(\bold{X})=\frac{X_2}{\alpha_2}$), as given in Table \ref{table:4}, Table \ref{table:5} and Table \ref{table:6}, receptively.
The following observations are evident from Tables 4-6:
\vspace*{2mm}

\noindent (i) All the GPN values are greater than $0.5$, confirming our theoretical findings of Example 4.2.1.
\\ (ii) When $\alpha_1$ is relatively larger than $\alpha_2$, we get higher GPN values.

%	\newpage
\begin{table}[h!]
	\centering
	\caption{The GPN values of the modified restricted MLE ($\delta_{2,RMLE}^{*}$) relative to the restricted MLE ($\delta_{2,RMLE}$)}
	\label{table:4}
	\centering
	\begin{tabular}{@{}ccccccc@{}}
		\toprule
		& \multicolumn{6}{c}{$(\alpha_1,\alpha_2)$}            \\
		$\theta_2 / \theta_1$ & (0.5,0.2) & (0.2,0.8) & (1,1) & (5,2) & (1,30) & (30,1) \\ \midrule
		1                   & 0.552     & 0.538     & 0.518 & 0.515 & 0.508  & 0.503  \\
		1.5                 & 0.613     & 0.567     & 0.603 & 0.640 & 0.519  & 0.736  \\
		2                   & 0.664     & 0.580     & 0.627 & 0.635 & 0.522  & 0.699  \\
		2.5                 & 0.692     & 0.592     & 0.634 & 0.605 & 0.522  & 0.666  \\
		3                   & 0.715     & 0.600     & 0.635 & 0.586 & 0.522  & 0.643  \\
		3.5                 & 0.731     & 0.595     & 0.627 & 0.568 & 0.519  & 0.622  \\
		4                   & 0.740     & 0.606     & 0.623 & 0.554 & 0.515  & 0.610 \\ \bottomrule
	\end{tabular}
	
\end{table}

% Please add the following required packages to your document preamble:
% \usepackage{booktabs}
\begin{table}[h!]
	\caption{The GPN values of the modified PNSEE ($\delta_{2,PNSEE}^{*}$) relative to the PNSEE ($\delta_{2,PNSEE}$)} 
	\label{table:5}
	\centering
	\begin{tabular}{@{}ccccccc@{}}
		\toprule
		& \multicolumn{6}{c}{($(\alpha_1,\alpha_2)$)}              \\
		$\theta_2 / \theta_1$ & (0.5,0.2) & (0.2,0.8) & (1,1) & (5,2) & (1,30) & (30,1) \\ \midrule
		1                     & 0.573     & 0.522   & 0.568 & 0.600 & 0.514  & 0.695   \\
		1.5                   & 0.612     & 0.543   & 0.595 & 0.631 & 0.522  & 0.686   \\
		2                     & 0.632     & 0.545   & 0.601 & 0.599 & 0.519  & 0.648   \\
		2.5                   & 0.645     & 0.548   & 0.596 & 0.575 & 0.514  & 0.620   \\
		3                     & 0.650     & 0.551   & 0.585 & 0.558 & 0.510  & 0.602   \\
		3.5                   & 0.659     & 0.549   & 0.581 & 0.544 & 0.508  & 0.591   \\
		4                     & 0.656     & 0.549   & 0.571 & 0.537 & 0.506  & 0.580  \\ \bottomrule
	\end{tabular}
\end{table}

\begin{table}[h!]
	\caption{The GPN values of the restricted MLE ($\delta_{2,RMLE}$) relative to the unbiased estimator ($\delta_{2,UE}$)}
	\label{table:6}
	\centering
	\begin{tabular}{@{}ccccccc@{}}
		\toprule
		& \multicolumn{6}{c}{$(\alpha_1,\alpha_2)$}              \\
		$\theta_2 / \theta_1$ & (0.5,0.2) & (0.2,0.8) & (1,1) & (5,2) & (1,30) & (30,1) \\ \midrule
		1                   & 0.702     & 0.569     & 0.628 & 0.648 & 0.523  & 0.758  \\
		1.5                 & 0.736     & 0.579     & 0.642 & 0.668 & 0.529  & 0.744  \\
		2                   & 0.745     & 0.579     & 0.641 & 0.627 & 0.522  & 0.698  \\
		2.5                 & 0.756     & 0.575     & 0.633 & 0.598 & 0.516  & 0.664  \\
		3                   & 0.751     & 0.578     & 0.616 & 0.575 & 0.512  & 0.640  \\
		3.5                 & 0.748     & 0.573     & 0.610 & 0.559 & 0.509  & 0.625  \\
		4                   & 0.744     & 0.572     & 0.597 & 0.550 & 0.506  & 0.611 \\   \bottomrule
	\end{tabular}
\end{table}

%	\section{\textbf{Concluding Remarks}}    \label{6}
%	For the problem of estimating order-restricted location/scale parameters of a general bivariate model, we used the GPN criterion to derive a better estimator than an arbitrary location/scale equivariant estimator. The GPN criterion is not yet used in the literature to obtain improved estimators for restricted parameters. We also unified various results by considering a general bivariate probability location/scale model (possibly having dependent marginals). Through the application of our findings, we obtained some new estimators (which perform better) for various specific probability models. We also generalized the results of various researchers, who compared various estimators using the Pitman nearer criterion, to the GPN criterion.

\section*{\textbf{Funding}}

This work was supported by the [Council of Scientific and Industrial Research (CSIR)] under Grant [number 09/092(0986)/2018].

\bibliography{Paper4}	
\end{document}